\documentclass[a4paper,leqno]{amsart}
\def\timestamp{July 11, 2001}
\usepackage{amsmath,amssymb,amscd,amsthm}

\date{\timestamp}
\makeatletter
\theoremstyle{plain}
 \newtheorem{theorem}{Theorem}[section]
 
 \newtheorem{proposition}[theorem]{Proposition}
 \newtheorem{lemma}[theorem]{Lemma}
 \newtheorem{corollary}[theorem]{Corollary}
 \newtheorem{assertion}[theorem]{Assertion}
 \newtheorem{fact}[theorem]{Fact}
\theoremstyle{definition}
 \newtheorem{definition}[theorem]{Definition}
 \newtheorem{remark}[theorem]{Remark}
 \newtheorem{example}[theorem]{Example}
\theoremstyle{plain}
\numberwithin{equation}{section}
\newcommand{\cmcone}{\mbox{\rm CMC}-$1$}
\newcommand{\ord}{\operatornamewithlimits{ord}}
\newcommand{\rank}{\operatorname{rank}}
\newcommand{\trace}{\operatorname{tr}}
\newcommand{\Res}{\operatornamewithlimits{Res}}
\newcommand{\dRes}{\displaystyle\operatornamewithlimits{Res}}
\newcommand{\diag}{\operatorname{diag}}
\newcommand{\bmath}[1]{\mbox{\boldmath $#1$}}
\newcommand{\smath}[1]{\mbox{\mathversion{bold}\scriptsize{$#1$}}}
\newcommand{\Z}{{\bmath{Z}}}
\newcommand{\R}{{\bmath{R}}}

\newcommand{\C}{{\bmath{C}}}
\newcommand{\CP}{{\bmath{C\!P}}}
\newcommand{\sR}{{\smath{R}}}
\newcommand{\GL}{\operatorname{GL}}

\newcommand{\SU}{\operatorname{SU}}

\newcommand{\SL}{\operatorname{SL}}
\newcommand{\Sl}{\operatorname{\mathfrak {sl}}}
\newcommand{\su}{\operatorname{\mathfrak {su}}}

\newcommand{\Ad}{\operatorname{Ad}}
\newcommand{\ad}{\operatorname{ad}}
\newcommand{\id}{\operatorname{id}}
\renewcommand{\Re}{\operatorname{Re}}
\renewcommand{\Im}{\operatorname{Im}}
\newcommand{\Hyp}{{\mathcal H}}
\newcommand{\Quad}{{\mathcal Q}}
\newcommand{\Lie}[1]{{\mathfrak {#1}}}
\newcommand{\TA}{\operatorname{TA}}
\newcommand{\trans}[1]{{\vphantom{#1}^t\!#1}}
\newcommand\al{a}
\newcommand\be{b}
\newcommand\zb{\bar{z}}
\newcommand\pd{\partial}
\newcommand\pdb{\overline{\partial}}
\makeatother
\title[An analogue of minimal surface theory]{
      An analogue of minimal surface theory \\
      in \boldmath$\SL(n,\C)/\SU(n)$ 
}
\author{M. Kokubu}
\address[Masatoshi Kokubu]{%
   Department of Natural Science,
   School of Engineering,
   Tokyo Denki University,
   2-2, Kanda-Nishiki-Cho,
   Chiyoda-Ku, Tokyo, 101-8457
   Japan
}
\email{kokubu@cck.dendai.ac.jp}
\author{M. Takahashi}
\address[Masaro Takahashi]{%
   Department of General Education,
   Kurume National College of Technology,
   Kurume, Fukuoka 830-8555, Japan
}
\email{taka@GES.kurume-nct.ac.jp}
\author{M. Umehara}
\address[Masaaki Umehara]{%
   Department of Mathematics, Faculty of Science,
   Hiroshima University,
   Higashi-Hiroshima 739-8526, Japan
}
\email{umehara@math.sci.hiroshima-u.ac.jp}
\author{K. Yamada}
\address[Kotaro Yamada]{%
   Faculty of Mathematics,
   Kyushu University 36,
   Hakozaki 6-10-1, Higashi-ku, Fukuoka 812-8581, Japan%
}
\email{kotaro@math.kyushu-u.ac.jp}
\subjclass{Primary 53A10; Secondary 53A35, 53A07}
\date{\timestamp}
\begin{document}
\begin{abstract}
  We shall discuss the class of surfaces with holomorphic
  right Gauss maps in non-compact duals of compact semi-simple Lie
  groups (e.g. \linebreak  $\SL(n,\C)/\SU(n)$), 
  which contains minimal surfaces in $\R^n$ and constant mean curvature
  $1$ surfaces in $H^3$.
  A Weierstrass type representation formula and
  a Chern-Osserman type inequality for such surfaces are given.
\end{abstract}
\maketitle
\section*{Introduction}
Surfaces of constant mean curvature one (\cmcone{})  in hyperbolic
$3$-space $\Hyp^3$
have the following three properties which are quite similar to properties
of minimal surfaces in $\R^3$:
\begin{enumerate}
\item The hyperbolic Gauss maps of \cmcone{}
      surfaces in $\Hyp^3$ are conformal, as is so for 
      the Gauss maps of minimal surfaces in $\R^3$.
\item Any simply connected regions of \cmcone{} surfaces in $\Hyp^3$ are 
      isometrically realized as minimal surfaces in $\R^3$.
\item An analogue of the Weierstrass representation formula
      (called the Bryant representation formula)
      for \cmcone{} surfaces in $\Hyp^3$ is known, 
      by which one can construct \cmcone{} surfaces via holomorphic data.
\end{enumerate}
Now \cmcone{} surface theory in $\Hyp^3$ is well developed
and many examples are known.
So it is quite natural to ask what
is the canonical generalization
of \cmcone{} surface theory for a much wider class of ambient spaces.
Hyperbolic $3$-space can be expressed as 
a quotient of a complex semi-simple Lie group  by its
compact real form, that is
\[
      \Hyp^3=\SL(2,\C)/\SU(2).
\]
However, higher dimensional hyperbolic spaces $\Hyp^n$ ($n\ge 4$)
have no such expressions. 
We shall show that a class of surfaces in certain non-compact type
symmetric spaces (containing $\SL(n,\C)/\SU(n)$) inherits all of the above
three properties, just as for
\cmcone{} surfaces in $\Hyp^3$, as follows. 
It is well-known that minimal surfaces in $\R^n$ 
can be lifted to null holomorphic maps
into the complex abelian group $\C^n$. 
Similarly,  our surfaces can be lifted to null holomorphic maps
into a complex semi-simple Lie group. 
So, roughly speaking, this gives a non-commutative version of Euclidean
minimal surface theory.

Let $G$ be a complex semi-simple Lie group and $H$ the compact real form
of $G$.
Then the quotient $G/H$ has the structure of a Riemannian symmetric space
and by the Cartan embedding we can regard $G/H\subset G$.
Via the adjoint representation,
we may assume that $G\subset \SL(n,\C)$ ($n=\dim G$).
Let
\[
  f\colon{}M\longrightarrow G/H (\subset G \subset \SL(n,\C))
\]
be a conformal immersion of a Riemann surface $M$ into $G/H$.
We define the right Gauss map
\[
  \nu_R\colon{}M\ni z \mapsto [f_zf^{-1}]\in P(\ad(\Lie{g}))
     (\subset P(\Sl(n,\C))\,),
\]
where $z$ is a complex coordinate and $f_z=\partial f/\partial z$.
Though $f_zf^{-1}\in \ad(\Lie{g})$ depends on the choice of
complex coordinate,  the projection $[f_zf^{-1}]$ into
the projective space $P(\ad(\Lie{g}))$ is globally defined on $M$.
We consider surfaces with holomorphic right Gauss maps.
By a direct calculation, it can be seen  that the square of the 
length of the mean curvature vector fields of such surfaces are
proportional to the sectional curvature of the ambient space
with respect to the tangent planes of the surfaces.
In particular, 
such surfaces in $\Hyp^3$
coincide with \cmcone{} surfaces.
We shall show the following:
\newtheorem*{introthm-A}{The canonical correspondence}
\begin{introthm-A}
 Any simply connected region of a
 minimal surface in $\R^n$ $(n=\dim G)$
 can be isometrically realized as a surface in $G/H$
 with holomorphic right Gauss map, and vice versa.
\end{introthm-A}

It should be remarked that minimal surfaces in $\R^r (\subset G/H)$ 
$(r:=\rank G)$ have holomorphic right Gauss maps whose 
 canonical correspondences are congruent to the original ones.
So the classical minimal surface theory is included in our theory. 

Let $\Lie{g}$ be the Lie algebra of $G$ and
$B$ the Killing form of $\Lie{g}$.

\newtheorem*{introthm-B}{Weierstrass-Bryant type representation formula}
\begin{introthm-B}
  Let $M$ be a simply connected Riemann surface and $\alpha$ 
  a $\Lie{g}$-valued holomorphic $1$-form such that 
\[
  B(\alpha,\alpha)=0,\qquad -B(\alpha, \sigma(\alpha))>0,
\]
  where $\sigma$ is the involution on $\Lie{g}$ with respect to
  the Riemannian symmetric pair $(G,H)$.
  Let  $F\colon{}M\to G$ be a solution of the ordinary
  differential equation
\[
     F^{-1}dF=\alpha.
\]
  Then the projection of $F$ into $G/H$
  gives a conformal immersion with holomorphic right Gauss map.
  Conversely, any conformal immersion with holomorphic right Gauss map
  is given in this manner.
\end{introthm-B}

Let $M$ be a Riemann surface and $f\colon{}M\to G/H$  a conformal
immersion with holomorphic right Gauss map. By the
above representation formula,  $f$ can
be expressed as the projection of the null holomorphic map
\[
  F\colon{}\widetilde{M}\longrightarrow G,
\]
where ``null'' means that the pull back of $B$ by $F$ vanishes, and
where $\widetilde M$ is the universal cover of $M$.
Then the inverse $F^{-1}\colon{}z \mapsto F(z)^{-1}$ of $F$ gives also a null holomorphic 
map and its projection
\[
  f^\#\colon{}\widetilde M\longrightarrow G/H
\]
gives a new conformal immersion with holomorphic right Gauss 
map. The immersion $f^\#$ is called the dual of $f$.
Though $f^\#$ is multi-valued on $M$, it has
single valued first fundamental form. 
If we now assume that the total absolute curvature of $f$
or $f^\#$ is finite, then there exists a compact Riemann
surface $\overline M$ and finite points $\{p_1,\dots,p_n\}$
such that $M$ is biholomorphic to $\overline M\setminus \{p_1,\dots,p_n\}$.
These points $\{p_1,\dots,p_n\}$ are called ends of the surface.
Similar to \cmcone{} surface theory in $H^3$, 
the total absolute curvature of $f^\#$ satisfies a much stronger
inequality than the Cohn-Vossen inequality:
 
\newtheorem*{introthm-C}{The Chern-Osserman type Inequality}
\begin{introthm-C}
 Let $f\colon{}M\to G/H$ be a conformal immersion with holomorphic right
 Gauss map whose first fundamental form is complete. 
 Suppose that the dual immersion $f^{\#}$ of $f$ has finite total
 absolute curvature $\TA(f^\#)<\infty$. 
 Then it satisfies the following inequality
\begin{equation}\tag{$*$}\label{eq:C-O}
   \frac{1}{2\pi}{\TA(f^\#)}\ge -\chi(M)+ r,
\end{equation}
 where $r$ is the number of ends
 and $\chi (M)$ is the Euler number of $M$.
\end{introthm-C}

For a \cmcone{} surface with finite total curvature 
in $\Hyp^3$, equality in \eqref{eq:C-O} implies that all ends are
regular and embedded (\cite{UY4}).
The condition for equality in \eqref{eq:C-O} in general case is still
unknown, whereas the equality condition of the Chern-Osserman
inequality for minimal surfaces in $\R^n$ is known (\cite{KUY}).

In the first section, we shall review the local properties of
minimal surfaces in $\R^n$.
In Section~\ref{sec:formula}, we shall describe the ambient space $G/H$
and give a fundamental formula for surfaces in it.
In Section~\ref{sec:gauss}, we shall show the canonical correspondence
and the Weierstrass-Bryant type representation formula.
In Section~\ref{sec:osserman}, we shall prove an Chern-Osserman type
inequality for surfaces with holomorphic right Gauss maps. 
For the proof, we apply the theory of complex ordinary differential
equation with regular singularities, which is prepared in Appendix. 
In Section~\ref{sec:perturb}, we shall construct a non-trivial 
example satisfying equality in the inequality \eqref{eq:C-O}
as a deformation of a minimal surface,
where techniques like in \cite{UY2} and \cite{ruy1} will
be applied.
\section{Preliminaries}
\label{sec:prelim}
In this section, we shall make a  quick review of the local properties
of minimal immersions of surfaces into $\R^n$, referring to several
facts from Lawson \cite{lawson}.

Let $M$ be a Riemann surface and $f=(f_1,\dots,f_n)\colon M \to \R^n$ 
a conformal minimal immersion, where $n\ge 3$ is an integer. 
For a complex coordinate $z$ on $M$,
we set $f_z:=(\partial f_j/\partial z)_{j=1,\dots,n}$
and $f_{\bar z}:=(\partial f_k/\partial \bar z)_{k=1,\dots,n}$, 
respectively.
We define a $\C^n$-valued 1-form  $\alpha$ on $M$ by
\[
    \alpha:=\partial f(=f_z dz), 
\]
which is independent of the choice of complex coordinates. 
Then the following two identities hold:
\begin{align}
   \overline{\partial}\partial f&= 0, 
        \label{eq:min}\\
   \langle\partial f, \partial f\rangle  &= 0
   \qquad(\text{that is } \langle f_z, f_z\rangle=0),  
        \label{eq:conf}
\end{align}
where $\langle~,~\rangle$ is the complexification of the canonical
inner product on $\R^n$.
The identity \eqref{eq:min} implies that the $\C^n$-valued $1$-form 
$\alpha$ is holomorphic. 

We set
\begin{align*}
    \alpha &= (\alpha_1,\dots,\alpha_n)\\
           &= \hat\alpha\, dz = (\hat\alpha_1,\dots,\hat\alpha_n)\,dz
\end{align*}
for a complex coordinate $z$ on $M$.
Then the induced metric $ds^2$ and the Gaussian curvature $K$ of
$f$ can be expressed in terms of $\alpha$ as follows:
\begin{align}
\label{eq:min-met}
   ds^2 &= \langle\alpha,\overline\alpha\rangle =
            \sum_{j=1}^n \alpha_j\cdot\overline{\alpha_j}, \\
\label{eq:min-gauss}
       K&= \left.\left(-\sum_{i<j}|\hat\alpha_i(\hat\alpha_j)_z
                        -\hat\alpha_j(\hat\alpha_i)_z|^2\right)
           \right/
           \left(\displaystyle\sum_{k=1}^n|\hat \alpha_j|^2\right)^3,
\end{align}
where the ``\,$\cdot$\,'' means the symmetric product.
In particular, $K\le 0$.

Let $\hat\pi\colon{}\widetilde M\to M$ be the universal covering
of the Riemann surface $M$ and fix the base point $z_0\in M$ and
$\tilde z_0\in\hat\pi^{-1}(z_0)$.
The holomorphic map 
\[
   F\colon{}\widetilde M\ni z\longmapsto
            \int_{\tilde z_0}^z\alpha\,\in\C^n
\]
is called a {\it lift\/} of $f$.
Then $f$ is written as $f=F+\overline F$.
In terms of the lift, condition \eqref{eq:conf} is equivalent to
\begin{equation}\label{null}
    \langle F_z, F_z\rangle=0.
\end{equation}
A holomorphic map $F$ of a Riemann surface into $\C^n$ is said to be
{\it null\/} if it satisfies condition \eqref{null}.

Conversely, suppose now that there is a $\C^n$-valued holomorphic
$1$-form on $M$ satisfying the following two conditions
\begin{align}
   \langle\alpha, \alpha\rangle &=0 
     \label{eq:conf2}, \\
   \langle\alpha , \overline{\alpha}\rangle &>0. 
     \label{eq:pos}
\end{align}
Then the map defined by
\begin{equation}\label{W}
   f:=2\Re\left(\int_{\tilde z_0}^z \alpha \right)\colon{}
                      \widetilde M\to \R^n
\end{equation}
is a conformal minimal immersion.
The condition \eqref{eq:conf2} (resp.~\eqref{eq:pos})
implies conformality (resp.~nondegeneracy of the induced metric)
of $f$.
The formula \eqref{W} is called  the {\it Weierstrass representation
formula}.

It should be remarked that
the immersion $f$ given by \eqref{W} might not be single-valued on $M$ in
general because the contour integral may depend on paths.
The immersion $f$ obtained by \eqref{W} is well-defined on $M$ if and
only if
\begin{equation}\label{period}
   \Re\int_{\gamma}\alpha=0\qquad 
       \text{for all loops $\gamma$ on $M$ emanating from $z_0$}.
\end{equation}

Finally, we remark on the {\it Gauss map\/} of minimal immersions:
We denote the complex hyperquadric by
\[
    \Quad_{n-2}:=
    \left\{[\xi]\in \CP^{n-1}\,;\, \langle\xi, \xi\rangle=0\right\}.
\]
Then the Gauss map of a conformal minimal immersion $f$ is defined by
\begin{equation}\label{eq:gauss-map}
     \nu\colon{}M\ni p\longmapsto 
       \left[ 
         \frac{\partial f_1}{\partial z}(p):
         \frac{\partial f_2}{\partial z}(p):\dots:
         \frac{\partial f_n}{\partial z}(p)\right]\in
         \Quad_{n-2},
\end{equation}
where the bracket denotes the homogeneous coordinate on $\CP^{n-1}$.
The map $\nu$ is independent of the choice of local
complex coordinate $z$. 
Let $d\sigma^2$ denote the Fubini-Study metric on $\CP^{n-1}$.
Then
\begin{equation}\label{eq:FS}
     (-K)ds^2=d\sigma^2
\end{equation}
holds.
\section{A fundamental formula  for surfaces in $\SL(n,\C)/\SU(n)$}
\label{sec:formula}
In analogy to minimal surface theory in $\R^n$, 
we shall construct a theory of surfaces in a certain class of 
Riemannian symmetric spaces of non-compact type $N$, which 
contains $\SL(n,\C)/\SU(n)$ as a typical example. 

\subsection*{The ambient symmetric space}
First we shall describe the ambient space $N$.
Let $\widetilde G$ be a simply-connected, complex semi-simple 
Lie group with Lie algebra $\Lie{g}$. 
We denote by $B$ the Killing form of $\Lie{g}$. 
Let $\Lie{h}$ denote a compact real form of $\Lie{g}$. 
Namely, $\Lie{h}$ is a real Lie subalgebra of $\Lie{g}$ 
satisfying
\begin{equation}\label{2:eq:decompo} 
  \Lie{g}=\Lie{h}  + \sqrt{-1}\,\Lie{h} 
\end{equation}
and the Killing form $B$ is negative definite on $\Lie{h}$.
Let $\sigma_0\colon \Lie{g}\to \Lie{g}$ be complex conjugation with
respect to the decomposition \eqref{2:eq:decompo}, which is an
involutive automorphism of $\Lie{g}$ as a {\it real\/} Lie algebra. 
Since $\widetilde G$ is simply connected, there exists a unique
involution 
\begin{equation}
  \tilde \sigma\colon{}\widetilde G \longrightarrow \widetilde G
\end{equation}
such that $(\tilde \sigma_*)_e$ coincides with $\sigma_0$,
where $(\tilde \sigma_*)_e$ is the differential of the map 
$\tilde \sigma$ at the identity $e \in \widetilde G$. 
We denote by $\widetilde H$ the set of all the fixed points of $\tilde\sigma$. 
It is known that $\widetilde H$ coincides with the unique connected 
Lie subgroup of $\widetilde G$ whose Lie algebra is $\Lie{h}$.
Also,
we define a positive definite inner product $g$ on $\sqrt{-1} \Lie{h}$ by 
\begin{equation}\label{2:eq:gandB}
    g(X,Y) = B(X, Y)
      \quad \text{for} \quad X, Y \in {\sqrt{-1}\,\Lie{h}}. 
\end{equation}
We let $\pi \colon G\to G/H$ be the canonical projection 
and $o=\pi(e)$.
It is well known that $(\widetilde G, \widetilde H; \tilde \sigma,g)$
is a Riemannian symmetric pair which defines a symmetric space of
non-compact type
\begin{equation}\label{eq:ambient}
    N=\widetilde G/\widetilde H
\end{equation}
in such a way that the tangent space $T_o N$ at the point $o$
is isomorphic to $\sqrt{-1}\Lie{h}$ as an inner product space. 
A typical example is $N=\SL(n,\C)/\SU(n)$. 
When $n=2$, $\SL(2,\C)/\SU(2)$ is isometric to the $3$-dimensional
hyperbolic space $\Hyp^3$. 

\subsection*{Canonical embedding of {\boldmath $N$} into 
             {\boldmath$\Ad(G)\,(\subset \SL(n,\C)$)}}

Consider the adjoint representation $\Ad\colon{}\widetilde
G\to\GL(\Lie{g})$. 
Since the Killing form $B$ is preserved under the adjoint action, 
the determinant of the linear transformation $\Ad(\tilde a)$ is $\pm 1$
for any $\tilde a \in \widetilde G$. 
Then 
\begin{equation}
  \det\left(\Ad(\tilde a)\right)=1 
\end{equation}
holds for all $\tilde a\in\widetilde G$,  because $\widetilde G$ is
connected.
Since $B$ is negative definite on $\Lie{h}$, we can take a basis
\begin{equation}\label{eq:basis-h}
   \{e_1,\dots,e_n\}\subset \Lie{h}\qquad \text{so that}\qquad
   B( e_i,e_j) =-\delta_{ij},
\end{equation}
where $n=\dim_{\sR}\Lie{h}$.
Then 
\begin{equation}\label{eq:basis}
   \{e_1,\dots,e_n\}\subset \Lie{g}
\end{equation}
is a basis of $\Lie{g}$ over $\C$.
With respect to this basis, we have 
\begin{equation}
  \Ad(\widetilde H)\subset \Ad(\widetilde G)\subset \SL(n,\C).
\end{equation}
If we put 
\[
    G:=\Ad(\widetilde G), \qquad H:=\Ad(\widetilde H), 
\]
then we get another expression 
\[
   N=G/H,\qquad G\subset \SL(n,\C), 
       \quad\text{and}\quad H\subset \SU(n). 
\]
This is verified using the fact that $\widetilde H$ contains the center
of $\widetilde G$ (\cite{helgason}, Theorem~1.1 in Chapter~VI).
This fact also implies that one can define an involution 
$\sigma \colon G\to G$ by 
\begin{equation}
  \sigma(\Ad(\tilde a))=\Ad(\tilde \sigma(\tilde a)) 
     \qquad (\tilde a\in \widetilde G). 
\end{equation}
Therefore we have another Riemannian symmetric pair $(G,H;\sigma,g)$ 
of $N$ such that $G$ is a Lie subgroup of $\SL(n,\C)$. 
This observation plays an important role for us, from the technical
viewpoint. 

Now we recall that $N= G/H$ can be embedded in $G$ as follows 
(\cite{helgason}, p.~276): 
The mapping $\psi \colon N=G/H \to G$ defined by 
\[
  \psi (aH) = a \sigma (a^{-1}) \quad \text{for} \quad a \in G, 
\]is a diffeomorphism of $N$ into $\psi(N)$. 
The image $\psi(N)$ coincides with the identity component 
of 
\[
    C_0:=\{a\in G\,;\, \sigma(a^{-1})=a\}.
\]
The mapping $\psi$ 
is often called the {\it Cartan embedding\/} and the image 
the {\it Cartan model}. 
\begin{lemma}\label{2:lem:sigma}
  For $X\in \Lie{g}$ $(\subset \Sl(n,\C))$ 
  and $a\in G(\subset \SL(n,\C))$, 
\[
   \sigma_0(X)=-X^*\quad\text{and}\quad \sigma(a^{-1})=a^*
\]
  hold, where ${}^*$ denotes the operation of conjugation and 
  transposition, namely $a^*=\trans{\bar a}$.
  Furthermore $\sigma_0(X)=-X$ holds for $X\in T_oN$.
\end{lemma}
\begin{proof}
  Write $X=X_1+\sqrt{-1}X_2$ ($X_1,X_2\in \Lie{h}$) with respect 
  to the decomposition \eqref{2:eq:decompo}. 
  Since $H\subset \SU(n)$, 
  $\Lie{h}$ is contained in the Lie algebra $\su(n)$ of $\SU(n)$. 
  Hence $(X_j)^*=-X_j$ ($j=1,2$). 
  This implies that 
\[
     \sigma_0 (X)=\sigma_0 (X_{1} + \sqrt{-1} X_{2})=
        X_{1} - \sqrt{-1} X_{2}=
       -X_{1}^{*} + \sqrt{-1} X_{2}^{*}=-X^{*}. 
\]
  The second assertion $\sigma(a^{-1}) = a^{*}$ follows from the
  well-known formula $\sigma (\exp X) = \exp \sigma_0 (X)$. 
  The last assertion follows from the fact $X=X^*$ for 
  $X\in T_oN=\sqrt{-1}\,\Lie{h}$.
\end{proof}

By Lemma~\ref{2:lem:sigma}, $\psi \colon N \to G$ is given by 
\[
   \psi(aH)=aa^*,
\]
and $\psi(N)$ is the identity component of 
$C_0=\{a \in G\,;\,a=a^* \}$. 

We will proceed by identifying $N$ with its image $\psi(N)$. 
Let $M$ be a Riemann surface.
We think of a map $f \colon M \to G/H$ as 
a map $f \colon M \to G$ with values in the identity component of
$C_0$. 
We give a left invariant metric $g_0$ on $G$ by
\begin{equation}\label{metric}
   g_{0}(X,Y) = -B(X, \sigma (Y))
  \quad \text{for} \quad X, Y \in \Lie{g} \cong T_eG.
\end{equation}
Then the metric $g$ in \eqref{2:eq:gandB} is related to $g_0$ by
\begin{equation}\label{pull}
  \psi^{*} g_{0} = 4 g.
\end{equation}
The Killing form $B_{\Sl}$ of $\Sl(n,\C)$ is given by
\begin{equation}\label{SL-killing}
   B_{\Sl}(X,Y)=2n \trace(XY)\qquad (X,Y\in \Sl(n,\C)).
\end{equation}
By \eqref{metric}--\eqref{SL-killing}
and  Lemma \ref{2:lem:sigma},
we have
\begin{align}\label{useful}
   g(X,Y) & =  \frac{1}{4}B(d\psi(X),d\psi(Y)) \\
   & =  \frac{1}{8n} B_{\Sl}(\ad(d\psi(X)),\ad(d\psi(Y)))
   \qquad (X,Y\in T_oN).\nonumber
\end{align}
The Lie group $G$ acts on $N=\psi(N)$ isometrically by
\begin{equation}\label{eq:isometry}
    v\longmapsto a v a^*\qquad (v\in \psi(N), a\in G)
\end{equation}
with respect to the metric $g$.

\subsection*{General formula for surfaces in \boldmath$N$}
Now we shall give a general representation formula for surfaces 
in $N$. 
Let $\widetilde M$ be a simply connected Riemann surface and
$f\colon \widetilde M \to N=G/H$ a smooth map. 
Then there exists a smooth lift $\varPhi\colon \widetilde M \to G$, 
i.e., $\varPhi$ is a smooth map satisfying $\varPhi \varPhi^*=f$ 
on $\widetilde M$. 
We remark that the lift $\varPhi$ is determined uniquely up to 
multiplying by a smooth $H$-valued map $h\colon \widetilde M
\to H$ on the right.
Let $\theta$ be the pull-back of the left-invariant 
Maurer-Cartan form of $G$ via $\varPhi$, namely, 
$\theta$ is the $\Lie{g}$-valued $1$-form $\varPhi^{-1} d \varPhi$ 
on $\widetilde M$. 
Denote by $\xi$ the $(1,0)$-part of $\theta$, and by $\eta$ the
$(0,1)$-part, i.e., 
\begin{equation}\label{eq:2:frenet}
\left\{
\begin{aligned}
  \xi  &= \theta^{1,0} = \varPhi^{-1} \pd \varPhi, \\
  \eta &= \theta^{0,1} = \varPhi^{-1} \pdb \varPhi. 
\end{aligned}
\right.
\end{equation}
Then we have 
\begin{lemma}\label{2:lem:conf}
  A smooth map $f\colon \widetilde M \to G/H$ is 
  conformal if and only if any lift $\varPhi$ of $f$ satisfies 
\begin{equation}\label{eq:conformal}
    B(\xi+\eta^*,\xi+\eta^*)=0. 
\end{equation}
  The condition \eqref{eq:conformal} does not depend 
  on the choice of a lift $\varPhi$. 
\end{lemma}
\begin{proof}
  Let $z=x+\sqrt{-1}y$ be a local coordinate on $\widetilde M$. 
  Then conformality is equivalent to
\begin{equation}\label{2:eq:fzfz}
  g(f_x,f_x)=g(f_y,f_y),\qquad g(f_x,f_y)=0 
  \quad \text{on}\quad T_{f(z)}N. 
\end{equation} 
  Since $G$ acts isometrically on $N$ as in \eqref{eq:isometry},
   \eqref{2:eq:fzfz} is rewritten as
\begin{align*}
  g(\varPhi^{-1}f_x(\varPhi^{-1})^*,\varPhi^{-1}f_x(\varPhi^{-1})^*)
   &=g(\varPhi^{-1}f_y(\varPhi^{-1})^*,\varPhi^{-1}f_y(\varPhi^{-1})^*), \\
  g(\varPhi^{-1}f_x(\varPhi^{-1})^*,\varPhi^{-1}f_y(\varPhi^{-1})^*)
   &=0. 
\end{align*} 
  By \eqref{useful}, it is equivalent to the condition 
  $B(\varPhi^{-1}f_z(\varPhi^{-1})^*,\varPhi^{-1}f_z(\varPhi^{-1})^*)=0$.
  On the other hand, by differentiating $f=\varPhi\varPhi^*$, we have
\begin{equation} \label{eq:fz}
  f_z=\varPhi_z\varPhi^*+\varPhi(\varPhi^*)_z
     =\varPhi_z\varPhi^*+\varPhi(\varPhi_{\bar z})^*.
\end{equation} 
  Hence condition \eqref{2:eq:fzfz} is written as
\begin{align*}
  0&=B(\varPhi^{-1}f_z(\varPhi^{-1})^*,\varPhi^{-1}f_z(\varPhi^{-1})^*)\\
   &=B(\varPhi^{-1}\varPhi_z + (\varPhi_{\zb})^*(\varPhi^{-1})^*,
     \varPhi^{-1}\varPhi_z + (\varPhi_{\zb})^*(\varPhi^{-1})^*),
\end{align*}
  which proves \eqref{eq:conformal}.

  Next we prove the latter part. 
  Let $\hat \varPhi$ be another lift on $\widetilde M$ and denote 
  $\hat \xi={\hat \varPhi}^{-1} \pd {\hat \varPhi}$, 
  $\hat \eta={\hat \varPhi}^{-1} \pdb {\hat \varPhi}$. 
  The lifts 
  $\hat \varPhi$ and $\varPhi$ are related by $\hat \varPhi = \varPhi h$ 
  for some function $h\colon \widetilde M \to H$. It follows that 
\begin{equation*}
   \hat \xi = h^{-1} \xi h + h^{-1} \pd h, \qquad 
   \hat \eta = h^{-1} \eta h + h^{-1} \pdb h. 
\end{equation*}
  Since $h^*=h^{-1}$, we have 
\begin{align*}
  {\hat \eta}^* 
    & = (h^{-1} \eta h + h^{-1} \pdb h)^* 
    = h^{-1} \eta^* h + (\pd h^{-1}) h 
    = h^{-1} \eta^* h + (h^{-1} \pd h h^{-1}) h \\
    & = h^{-1} \eta^* h - h^{-1} \pd h. 
\end{align*}
  Therefore 
  $\hat \xi +{\hat \eta}^*=h^{-1}(\xi + \eta^*)h$. 
  This implies that 
\begin{equation*}
  B(\hat \xi+\hat \eta^*,\hat \xi+\hat \eta^*)=
  B(\xi+\eta^*,\xi+\eta^*).\qquad \qedsymbol
\end{equation*}
\renewcommand{\qedsymbol}{\relax}
\end{proof}
\begin{lemma}\label{2:lem-add}
  Let $f\colon \widetilde M \to G/H$ be a conformal map. 
  Then the induced metric  $ds^2$ is written 
  in terms of a lift $\varPhi$ as 
\begin{equation}\label{eq:indmet}
    ds^2 = B(\xi+\eta^*,\ \xi^*+\eta).
\end{equation}
  Moreover, $f$ is an immersion if and only if 
\begin{equation}\label{eq:immersion}
  B(\xi+\eta^*,\ \xi^*+\eta)>0.
\end{equation}
\end{lemma}
\begin{proof}
  By the conformality, the induced metric is
  $ds^2=4g(f_z,f_{\zb})dzd\zb$. 
  The proof is similar to that of Lemma \ref{2:lem:conf}. 
\end{proof}

The integrability condition of the differential equation 
\eqref{eq:2:frenet} is given by 
\begin{equation}\label{eq:intcond}
  \pdb \xi + \pd \eta + \eta \wedge \xi 
         + \xi \wedge \eta =0, 
\end{equation}
which is equivalent to 
\begin{equation}\label{eq:intcond2}
d\theta+\theta \wedge \theta=0. \tag{$\ref{eq:intcond}'$}
\end{equation}
So the following fundamental theorem holds:
\begin{theorem}\label{Thm_G}
  Let $\theta$ be a $\Lie{g}$-valued $1$-form on a 
  simply-connected Riemann surface $\widetilde M$ which satisfies 
\begin{equation}\label{G-C}
   d\theta+\theta \wedge \theta=0 
\end{equation}
  and 
\begin{equation*}
B((\theta+\theta^*)^{1,0},\ (\theta+\theta^*)^{1,0})=0. 
\end{equation*}
  Then there exists a conformal map $f\colon \widetilde M \to G/H$ with 
  the induced metric 
\begin{equation*}
  ds^2=B((\theta+\theta^*)^{1,0},\ (\theta+\theta^*)^{0,1}),
\end{equation*}
  which is determined up to $G$-congruence.  
\end{theorem}
\begin{proof}
By the assumption of integrability, there exists a solution 
$\varPhi \colon \widetilde M \to G$ of the differential equation 
$\varPhi^{-1}d \varPhi = \theta$ for an initial value $\varPhi(z_0)\in G$. 
Since $(\theta+\theta^*)^{1,0}=\xi+\eta^*$
and  $(\theta+\theta^*)^{0,1}=\xi^*+\eta$,
the mapping $f:=\varPhi \varPhi^*$ 
has the first fundamental form 
$B((\theta+\theta^*)^{1,0},(\theta+\theta^*)^{0,1})$,
by Lemma~\ref{2:lem-add},
and thus it is a conformal immersion.
Let $\hat \varPhi$ be a solution of the differential equation  
$\varPhi^{-1}d \varPhi = \theta$ under another initial condition. 
Then $\hat \varPhi$ differs from $\varPhi$ only by 
$\hat \varPhi = a \varPhi $ for some $a \in G$. 
Therefore, $\hat f:= \hat \varPhi \hat \varPhi^*$ satisfies 
$\hat f = afa^*$, which implies that $\hat f$ is congruent to $f$ 
by the isometric action of $a \in G$.  
\end{proof} 
\begin{remark}
When $N=\SL(2,\C)/\SU(2)$, that is, $N$ is the hyperbolic $3$-space,
the above theorem is interpreted as a generalization of
the fundamental theorem for surfaces:
If we choose $\varPhi$ as the Frenet frame
of the surface, the formula \eqref{eq:2:frenet} is exactly the Frenet
formula for the surface and the integrability condition
\eqref{G-C} is equivalent to the Gauss and Codazzi equations. 
The above formula also involves the Bryant representation formula
for \cmcone{} surfaces in $\Hyp^3$.
In fact, if we take $\theta$ as a holomorphic $1$-form
on $\widetilde M$, then the corresponding surface $f$ has
holomorphic Gauss map (See Theorem~\ref{Thm-B} in the next section)
and thus it has constant mean curvature $1$, and this is 
called the Bryant representation formula.
\end{remark}
\section{Surfaces with holomorphic right Gauss map}
\label{sec:gauss}
In this section, we shall define a ``right Gauss map'' for surfaces in
$N$ and discuss fundamental properties of surfaces with holomorphic
right Gauss maps.
Since \cmcone{} surfaces in $\Hyp^3$ can be characterized by
the holomorphicity of the Gauss map, it is a generalization of 
\cmcone{} surfaces in $\Hyp^3$.
We will use the same notation as in the previous section. 
\subsection*{Right Gauss maps}
Let $N=G/H$ be as in the previous section.
The Lie algebra $\ad(\Lie{g})$ of $G$ is a subalgebra of $\Sl(n,\C)$.
We identify $N$ with the image $\psi(N)\subset G\subset\SL(n,\C)$ of the
Cartan embedding $\psi$.
We denote by $P(\ad(\Lie{g}))$ the projective space of the complex
vector space $\ad(\Lie{g})$.

Let $f\colon{}M\to N=G/H$ be a conformal immersion.
Then for a complex coordinate $z$ on $M$, 
$f_z f^{-1}\in \ad(\Lie(g))$ for each point $z\in M$,
and here the product $f_z$ and $f^{-1}$ is matrix multiplication in
$\Sl(n,\C)$.
Although this map does depend on the choice of a coordinate $z$, the 
projection $[f_z f^{-1}]$ is independent of the choice of coordinate.
\begin{definition}
  For a conformal immersion $f\colon{}M\to G/H$,
  we define
\[
    \nu_R\colon{} M \ni z\longmapsto
                    [f_z f^{-1}]\in P(\ad(\Lie{g}))\subset
                                    P(\Sl(n,\C)).
\]  
  We call $\nu_R$ the {\it right Gauss map\/} of $f$.
\end{definition}
\begin{remark}\label{rem:gauss-quad}
  Let
\[
   \Quad_{n-2}(\ad(\Lie g)):=
      \{ [\xi]\in P(\ad(\Lie g))
                      \,; B(\xi,\xi)=0\},
\]
  that is, $\Quad_{n-2}(\ad(\Lie{g}))$ is the hyperquadric in
  $P(\ad(\Lie{g}))$  with respect to $B$.
  If $f\colon{}M\to G/H$ is a conformal immersion, $g(f_z,f_z)=0$
  holds. Hence the right Gauss map $\nu_R$ has values in $\Quad_{n-2}$.

  Using the basis $\{e_1,\dots,e_n\}$ of $\Lie{g}$ as in
  \eqref{eq:basis}, $\Quad_{n-2}(\ad(\Lie g))$ is identified with the 
  hyperquadric $\Quad_{n-2}$ in $\CP^{n-1}$.
  Thus, the right Gauss map can be considered as corresponding to the
  Gauss map of a minimal surface in $\R^n$ defined in
  \eqref{eq:gauss-map}.
\end{remark}

\subsection*{Bryant type representation formula}
Let $F\colon M \to G$ be a holomorphic map. 
We denote by $\alpha_F$ the pull-back of the left-invariant 
Maurer-Cartan form via $F$ and by $\alpha_F^\#$ the pull-back 
of the right-invariant Maurer-Cartan form via $F$, alternatively
\[
   \alpha_F = F^{-1}dF,
   \qquad
   \text{and}\qquad
   \alpha_F^{\#} = dFF^{-1}
\]
in matrix form. 
\begin{definition}
  A holomorphic map $F \colon M \to G$ is said to be {\it null\/} 
  (or {\it isotropic}) if 
  the pull back of the Killing form $B$ by $F$ vanishes, that is, 
\[
     B(\alpha_F,\alpha_F) \equiv 0.
\]
This is equivalent to the condition  
\[
     B(\alpha^\#_F,\alpha^\#_F) \equiv 0.
\]
\end{definition}
Let $\pi\colon{}G\to N=G/H$ be the canonical projection.
\begin{theorem}[Bryant type formula]\label{Thm-B} 
  Let $\widetilde M$ be a simply-connected Riemann surface
  and $\alpha$ be a $\Lie{g}$-valued holomorphic $1$-form
  on $\widetilde M$ satisfying the following two properties{\rm :}
\begin{enumerate}
  \item $B(\alpha,\alpha)$ vanishes identically on $\widetilde M$,
  \item $B(\alpha,\alpha^*)$ is positive definite on $\widetilde M$.
\end{enumerate}
  Let $F \colon \widetilde M\to G$ be a solution of the ordinary differential
  equation 
\begin{equation}\label{eq:repr-ode}
  F^{-1}dF=\alpha.
\end{equation} 
  Then $f=\pi\circ F \colon \widetilde M\to N$
  is a conformal immersion which has holomorphic right Gauss map.
  Conversely, any conformal immersion of $\widetilde M$ into $N$ with
  holomorphic right Gauss map is constructed in this manner.
\end{theorem}
\begin{proof}
  We set $\theta=\alpha$, so 
  it satisfies the integrability condition \eqref{G-C}.
  Thus by Theorem~\ref{Thm_G}, there exists a conformal 
  immersion $f \colon \widetilde M\to N$ such that
  $f=\pi\circ F =FF^*$, where $F$ is a solution of 
  the ordinary differential equation of $F^{-1}dF=\alpha$. 
  Since $\alpha$ is holomorphic, $F$ is also holomorphic.
  Now we have 
\[
    f_zf^{-1}=(FF^*)_z(FF^*)^{-1}=F_zF^{-1}=\hat \alpha
    \qquad (\alpha=\hat\alpha \,dz),
\]
  which implies the Gauss map $\nu_R$ of $f$ is holomorphic.

  Next we shall prove the converse.
  Let $f \colon \widetilde M\to N$ be a conformal immersion with
  holomorphic right Gauss map.
  Let us choose a lift $\varPhi \colon \widetilde M\to G$ 
  of $f$ arbitrarily.
  We set
\begin{equation}\label{eq:frenet}
\left\{
\begin{aligned}
  \xi  &= \varPhi^{-1}\pd\varPhi \\
  \eta &= \varPhi^{-1}\pdb\varPhi,
\end{aligned}
\right.
\end{equation}
  and assume $\xi$ and $\eta$ has a local expression
  $\xi=\hat\xi dz$, $\eta=\hat \eta dz$ for a local complex
  coordinate $z$. The integrability condition \eqref{eq:intcond}
  is written as 
\begin{equation}\label{3:eq:intcond}
  \hat\xi_{\zb}-\hat\eta_z+\hat \eta \hat \xi -\hat \xi \hat\eta =0,  
\end{equation}
  the conformality \eqref{eq:conformal} is  
\begin{equation}\label{3:eq:conformal}
    B(\hat\xi+\hat\eta^*,\hat\xi+\hat\eta^*)=0
\end{equation}
  and the condition of nondegeneracy \eqref{eq:immersion} is 
\begin{equation}\label{3:eq:immersion}
    B(\hat\xi+\hat\eta^*,\hat\xi^*+\hat\eta) > 0. 
\end{equation}

  On the other hand, we have 
\begin{align*}
   f_z f^{-1} &= (\varPhi \varPhi^*)_z(\varPhi \varPhi^*)^{-1}\\
              &= \varPhi_z \varPhi^{-1} + 
                 \varPhi (\varPhi^{-1} \varPhi_{\bar z})^* \varPhi^{-1}
                = \varPhi (\hat \xi + \hat\eta^*)  \varPhi^{-1}.
\end{align*}
  Hence the right Gauss map $\nu_R=[f_zf^{-1}]$ is holomorphic if and only
  if there exists a holomorphic map $\phi\colon{}\widetilde M\to\Lie g$
  and a non-vanishing function $\lambda\colon{}U\to\C^*$ such that
\begin{equation}\label{3:eq:ghol}
   \varPhi(\hat\xi+\hat\eta^*)\varPhi^{-1} = \lambda\phi.
\end{equation}
Differentiating \eqref{3:eq:ghol} with respect to $\zb$, 
we have 
\begin{equation*}
  \hat\eta\hat\xi + \hat\eta\hat\eta^* + \hat\xi_{\zb} + 
    (\hat\eta_z)^* -\hat\xi\hat\eta -\hat\eta^*\hat\eta = 
    (\log \lambda)_{\zb}(\hat\xi +\hat\eta^*). 
\end{equation*}
It follows from \eqref{3:eq:intcond} that 
\begin{equation*}
  \hat\eta_z + (\hat\eta_z)^* +\hat\eta\hat\eta^* -\hat\eta^*\hat\eta = 
     (\log \lambda)_{\zb}(\hat\xi+\hat\eta^*). 
\end{equation*}
Since the left-hand side of the above equation is Hermitian, 
so is the right-hand side, 
namely, 
\begin{equation*}
  (\log \lambda)_{\zb}(\hat\xi+\hat\eta^*)=
     \overline{(\log \lambda)_{\zb}}(\hat\xi^*+\hat\eta).
\end{equation*}
Hence, 
\begin{equation*}
  B\left((\log \lambda)_{\zb}(\hat\xi+\hat\eta^*),\hat\xi^*+\hat\eta\right)=
  B\left(\overline{(\log \lambda)_{\zb}}
       (\hat\xi^*+\hat\eta),\hat\xi^*+\hat\eta\right).
\end{equation*}
  This implies from \eqref{3:eq:conformal} and \eqref{3:eq:immersion} that 
  $(\log \lambda)_{\zb}=0$, i.e., 
  $\lambda$ must be holomorphic.
Therefore
\begin{equation}\label{eq:alpha-adj}
    \alpha^{\#}:=\varPhi(\xi+\eta^*)\varPhi^{-1}=\lambda\phi\,dz
\end{equation}
  is a holomorphic $1$-form on $\widetilde M$.
  Moreover, by conformality \eqref{3:eq:conformal}, $\alpha^{\#}$ is null.

  Next, let $F$ be a holomorphic map obtained by solving 
  the differential equation 
\begin{equation}\label{eq:ode-dual}
    dF F^{-1}=\alpha^{\#}
\end{equation}
  for $\alpha^{\#}$  as in \eqref{eq:alpha-adj}
  under the initial condition $F(z_0)=\varPhi(z_0)\in G$ for a base
  point $z_0\in \widetilde M$.
  We want to show that $F$ is also a lift of $f$, that is,
  $f=FF^{*}=\varPhi\varPhi^*$. 
  For this, we have only to prove that $F$ differs from $\varPhi$ by 
  an $H$-valued function, namely, that $\varPhi^{-1} F$ takes values 
  in $H$. 
  In fact, $\varPhi^{-1}F$ satisfies 
\begin{align*}
   d(\varPhi^{-1} F)(\varPhi^{-1} F)^{-1} 
     &= \{ -\varPhi^{-1}d\varPhi \varPhi^{-1} F 
        + \varPhi^{-1} dF \} (\varPhi^{-1} F)^{-1}
      \\
     &=-(\xi +\eta )+(\xi +\eta^*) = -\eta+\eta^*
\end{align*}
  and the right-hand side $-\eta+\eta^*$ is a $1$-form which has values 
  in $\Lie{g} \cap \su(n)$, i.e. $\Lie{h}$. 
  Therefore $\varPhi^{-1} F$ remains in $H$. 
  Thus $f$ admits a holomorphic lift $F$ if $\nu_R$ is holomorphic. 
\end{proof}
\begin{remark}\label{rem:ambiguity}
  For given conformal immersion $f\colon{}\widetilde M\to N$ with 
  holomorphic right Gauss map, the holomorphic immersion
  $F\colon{}\widetilde M \to G$ as in Theorem~\ref{Thm-B} is called the
  {\em holomorphic lift\/} of $f$.
  The holomorhic lift of $f$ is determined up to right multiplication
  by constant elements in $H$.
  In fact, let $F$ and $\widetilde F$ be two holomorphic lifts of $f$.
  Then $b:=F^{-1}\widetilde F$ is a map into $H$, because
  $FF^{*}=\widetilde F\widetilde F^{*}$.
  Holomorphicity of $F$ and $\widetilde F$ implies $b$ is a holomorphic
  map.
  On the other hand,  since $b\in H\subset\SU(n)$, $b^*=-b$ is also
  holomorphic. Hence $b$ is anti-holomorphic.
  This implies that $b\colon{}\widetilde M\to H$ is a constant map,
  and $\widetilde F=Fb$.
\end{remark}

As an application of the representation formula, we prove 
\begin{proposition}
  There are no compact surfaces without boundary
  with holomorphic right Gauss map in $N$. 
\end{proposition}
\begin{proof}
  Suppose that $M$ is a compact Riemann surface and
  $f\colon M \to G/H$ a conformal immersion with holomorphic 
  right Gauss map.
  Take  a holomorphic lift $F$ of $f$.

  We may assume that $G$ is a subgroup of $\SL(n,\C)$, and $N=G/H$
  is a subset in $\SL(n,\C)$, by Cartan embedding.
  Then the trace of $f\in \SL(n,\C)$ satisfies 
\[
   (\trace f)_{z\zb}=\trace (f_{z\zb})=\trace\left\{F_z(F_z)^*\right\}  
      \ge 0, 
\]
  where $z$ is a complex coordinate of $M$.
  Hence the function $\trace f\colon{}M\to\R$ is subharmonic,
  which must be constant since $M$ is compact. 
  By an isometry in $N$, we may assume that $f(z_0)=e$, where 
  $e$ is the identity element of $G$. 
  Then $\trace f$ is identically $n$. 
  On the other hand, $\det f$ is identically $1$. 
  This implies that 
  the eigenvalues $\lambda_1, \dots, \lambda_n$ of $f$ satisfy  
\[
    \lambda_1+ \dots + \lambda_n =n, \qquad \lambda_1 \cdots \lambda_n =1. 
\]
  Here, $\lambda_1,\dots,\lambda_n\in\R$, because $f$ is Hermitian.
  Moreover, they are positive around $z_0$ because all of them are $1$ 
  at $z_0$. 
  Therefore, 
\begin{equation}\label{eq:eigen} 
   \frac{\lambda_1+ \dots + \lambda_n}{n} 
        \ge \sqrt[n]{\lambda_1 \cdots \lambda_n} 
\end{equation}
  holds. However, since both sides are $1$, the equality is 
  attained in \eqref{eq:eigen}. 
  Thus, 
\[
      \lambda_1 = \cdots = \lambda_n =1. 
\]
  Since $f$ is Hermitian, this implies that $f(z)$ is equal to the
  identity matrix, a contradiction.
\end{proof}

\subsection*{A duality on surfaces with holomorphic right Gauss maps}
Let $f \colon M\to G/H$ be a conformal
immersion of a Riemann surface $M$ with holomorphic right Gauss map
and $F \colon \widetilde M\to G$ a holomorphic lift of $f$.
Then the inverse $F^{-1} \colon z \mapsto F(z)^{-1}$ is also
a null holomorphic immersion.
In particular,
\[
   f^{\#}:=\pi\circ F^{-1} \colon \widetilde M \longrightarrow G/H
\]
is a conformal immersion defined on $\widetilde M$
which has holomorphic right Gauss map.
The map $f^\#$ is called a {\it dual\/} of $f$.

The notion of duality is also defined for  minimal immersions
in $\R^n$. But the dual of $f$ is nothing but the antipodal
immersion $-f$ because the inverse of the lift $F$ is $-F$.
The dual of surfaces of holomorphic right Gauss map is not 
congruent to the original surface because the group $G$ is
non-commutative.

Since the holomorphic lift $F$ has an ambiguity of the right action of
$H$, the dual $f^{\#}$ depends on the choice of $F$.
If $\hat F=Fb$ ($b\in H$) is another choice of
a holomorphic lift, $\hat f^\#=\pi\circ \hat F^\#=bf^*b^*$ is congruent to
$f^\#$.
However, the duals of two congruent immersions $f$ and $\tilde f$ might not
be congruent. In fact, let $\tilde f:=afa^*$ ($a\in G\subset\SL(n,\C)$),
which  is congruent to $f$, by \eqref{eq:isometry}.
Then one can take the lift $\tilde F=a F$ of $\tilde f$. Hence
$\tilde f^{\#}=F^{-1}a^{-1}(a^{*})^{-1}F^*$ is not congruent to $f$ 
in general.

Let $M$ be a (not necessarily simply-connected) Riemann surface 
and $f\colon{}M\to N$ a conformal immersion with holomorphic right 
Gauss map.
Then a $\Lie{g}$-valued $1$-form $\alpha_F$ is only defined on 
the universal cover $\widetilde M$ of $M$
(see Example~\ref{3:ex:cate}), even though $\alpha$ for a
conformal minimal immersion $f\colon{}M\to\R^n$ is well-defined on $M$.
On the contrary, $\alpha_F^{\#}$ is well-defined on $M$, as seen in the
following:
\begin{proposition}
  Let $f{}\colon M\to G/H$ be a conformal immersion
  with holomorphic right Gauss map,
  and let $F \colon \widetilde M\to G$ be its holomorphic lift.
  Then it holds that
\[
    \alpha^\#_F=-\alpha_{F^{-1}}.
\]
  Moreover, the $\Lie{g}$-valued null holomorphic $1$-form $\alpha^\#_F$
  is single-valued on $M$.
\end{proposition}
\begin{proof}
  In fact, we have
\[
  \alpha_{F^{-1}}=F d(F^{-1})=F(-F^{-1}dF F^{-1})=-dFF^{-1}=-\alpha^\#_F, 
\]
  which proves the first assertion. 
  Let $\tau$ be a covering transformation of the universal covering 
  $\widetilde M$ of $M$. 
  Since $F \circ \tau$ also induces $f$, there exists $b_{\tau}\in H$
  such that 
  $F \circ \tau(z) =F(z)b_{\tau}$ (see Remark~\ref{rem:ambiguity}). 
  Hence
\[
  \alpha^\#_F\circ \tau=d(F\circ \tau) (F\circ \tau)^{-1}
     =dF b_{\tau} (F b_{\tau})^{-1}=dFF^{-1}=\alpha^\#_F.
\]
  This shows that $\alpha^\#_F$ is single-valued on $M$.
\end{proof}

We shall use the notation ${}^\#$ for the differential geometric
invariants induced by $f^\#$ throughout this paper. 
\begin{corollary}\label{prop:metric}
  Let $f \colon M\to G/H$ be a conformal  immersion with
  holomorphic Gauss map, and $f^\#$ its dual. 
  Then $f^\#$ induces a well-defined metric 
  ${ds^2}^{\#}$ on $M$. 
\end{corollary}
\begin{proof}
By \eqref{metric2}, we have
\begin{equation}\label{metric3}
  {ds^2}^\#=\sum_{j=1}^n \alpha_j^\#\cdot \overline{\alpha_j^\#}.
\end{equation}
  Since $\alpha_F^\#$ is single-valued, this implies ${ds^2}^\#$
  is single-valued.
\end{proof}

\subsection*{The canonical correspondence}
Let $f\colon{}\widetilde M\to N$ be a conformal immersion with holomorphic 
right Gauss map and let $F$ be a holomorphic lift of it.
Take a basis $\{e_1,\dots,e_n\}$ of $\Lie{h}$ as in \eqref{eq:basis-h}.
Then $\alpha_F=F^{-1}dF$ is written as
$\alpha_F=\sum \alpha_je_j$ for some holomorphic $1$-forms $\alpha_j$ on 
$\widetilde{M}$ ($j=1,\dots,n$).  
One can easily verify that $F$ is null if and only if  
\begin{equation}\label{null2}
  \sum_{j=1}^n \alpha_j\cdot \alpha_j=0,
\end{equation}
and the induced metric $ds^2$ of $f$ is written as
\begin{equation}\label{metric2}
  ds^2=\sum_{j=1}^n \alpha_j\cdot \overline{\alpha_j}.
\end{equation}
Now we define a map 
\begin{equation}\label{can-corr}
  f_0:=2\Re\int_{z_0}^z (\alpha_1,\dots,\alpha_n )
            \colon \widetilde M\longrightarrow \R^n.
\end{equation}
By \eqref{null2} and \eqref{metric2}, $f_0$ 
is a conformal minimal immersion on the universal covering
$\widetilde M$, which is locally isometric to $f$. 
The minimal immersion $f_0$ is called the {\it canonical
correspondence\/} of $f$.
Since $f$ and $f_0$ have the same first fundamental form,
local intrinsic properties of the two surfaces are the same.
For example, by \eqref{eq:min-gauss}, the Gaussian curvature $K$ of
$f\colon{}M\to N$ is a non-positive real function on $M$.

On the other hand, 
the canonical correspondence $f^{\#}_0$ of the dual $f^{\#}$
has the same Gauss map as $f$ in the sense of
Remark~\ref{rem:gauss-quad}, but is not isometric to $f$ in general.

As an application of the canonical correspondence, 
we shall prove the following fact.
\begin{proposition}\label{prop:cp-fi}
  Let $f \colon M \to G/H$ be a complete conformal immersion with
  holomorphic right Gauss map of finite total curvature or dual finite 
  total curvature, that is,
\[
    \int_M (-K)\,dA < + \infty \qquad\text{or}\qquad
    \int_M (-K^{\#})\,dA^{\#}<+\infty,
\]
  where $K$ {\rm(}resp.~$K^{\#}${\rm)} is the Gaussian curvature
  and $dA$ {\rm(}resp.~$dA^{\#}${\rm)} is the area element with respect 
  to the induced metric $ds^2$ {\rm(}resp. $ds^2{}^{\#}${\rm)}.
  Then $M$ is biholomorphic to a compact Riemann surface with
  finitely many points removed.
\end{proposition}
\begin{proof}
  Suppose that $f$ is of finite total curvature (resp.~finite
  dual total curvature). Since the metric induced by $f$ (resp.~$f^\#$)
  is locally isometric to the metric of a minimal 
  surface in Euclidean space, it has non-positive 
  Gaussian curvature on $M$. 
  So $M$ admits a complete metric of non-positive curvature,
  which yields the assertion by the same argument as 
  in minimal surface theory. (See \cite{lawson}.)
\end{proof}
Under the assumption of Proposition~\ref{prop:cp-fi}, we 
may write $M =\overline M \setminus \{p_1,\dots,p_r\}$ 
where $\overline M$ is a compact Riemann surface. 
Each $p_j$ is called an end. 
By definition, the Euler number $\chi(M)$ of $M$ 
is equal to $\chi(\overline M)-r$.  

\subsection*{Examples}
Finally in this section, we give simple examples of a complete surfaces
with holomorphic right Gauss maps and  with finite total absolute
curvature.

\begin{example}\label{ex-add1}
We denote by $r$ the rank of the Lie group $G$.
Then the Euclidean space $\R^r$ is a totally geodesic submanifold in
 $G/H$.
Let $M$ be a Riemann surface and $f_0\colon{}M^2\to \R^r$
a conformal minimal immersion.
Then $f_0$ has holomorphic right Gauss map as a surface in $G/H$. 
Here we demonstrate it for $G=\SL(r+1,\C)$.
Let $D$ be the maximal abelian subgroup in 
$\SL(r+1,\C)$ which consists of diagonal matrices in $\SL(r+1,\C)$.
Then the projection $\pi(D)$ of $D$ into $\SL(r+1,\C)/\SU(r+1)$
is isometric to the Euclidean space $\R^r$.
Let $f\colon{}M^2 \to D$ be a conformal immersion of the form 
\[
  f=\diag(f_1,...,f_{r+1}),
\]
where $\diag(f_1,\dots,f_{r+1})$ denotes the diagonal
matrix whose diagonal components are $f_1,\dots,f_{r+1}$.
Then we have
\[
   f^{-1}df_=\diag(d\log f_1,...,d\log f_{r+1}).
\]
Thus $f^{-1}df$ is holomorphic if and only if
$(d\log f_1,...,d\log f_{r+1})$ is also, which implies
$f$ is a conformal minimal immersion in $D$
if and only if it has holomorphic right Gauss map 
as a surface in $\SL(r+1,\C)/\SU(r+1)$.
\end{example}

\begin{example}\label{ex-add2}
Let $M$ be a Riemann surface and $a,b,c$ meromorphic functions
on $M$. 
We set
\[
  F:=\begin{pmatrix}
     1 & a &b \\
     0 & 1 & c \\
     0 & 0 & 1 
     \end{pmatrix}.
\]
Then $F$ is a null meromorphic map on $M^2$.
So we can construct many surfaces having holomorphic right
Gauss maps as the projection of such $F$.
For example, if we set $M=\C\setminus\{0\}$ and
$a(z)=1/b(z)=z$, $c(z)=1$, then it gives a complete
surface with finite total curvature in $\SL(3,\C)/\SU(3)$.
For a suitable choice of orthonormal basis,
one can easily check that the canonical correspondence 
$F_0\colon{}\widetilde M^2\to \C^8$ is given by
\[
   F_0:=\int_{z_0}^z
          (db-adc,i(db-adc),da,-ida,dc,idc,0,0).
\]
In particular, $F_0$ is not single-valued on $M^2$.
This implies that the canonical 
correspondence \eqref{can-corr} is local in nature.
\end{example}

\begin{example}\label{3:ex:cate}
  Let $M=\C\setminus\{0\}$, and let $\widetilde M$ be the universal cover 
  of $M$.
  Define a holomorphic map 
  $F_{\mu,\al,\be} \colon \widetilde{M} \to \SL(3,\C)$ by
\begin{equation}\label{3:eq:Fmuab}
  F_{\mu,\al,\be}(z)=
  \begin{bmatrix} 
    \sqrt{\frac{\be^2+3\mu^2}{\be^2-\al^2}} z^{\mu+\al} & 0 & 
    \sqrt{\frac{\al^2+3\mu^2}{\be^2-\al^2}} z^{\mu-\be} \\
                            0   & z^{-2\mu} & 0 \\
    \sqrt{\frac{\al^2+3\mu^2}{\be^2-\al^2}} z^{\mu+\be} & 0 & 
    \sqrt{\frac{\be^2+3\mu^2}{\be^2-\al^2}} z^{\mu-\al} 
  \end{bmatrix},
\end{equation}
where $\mu$, $\al$ and  $\be$ are real constants such that 
$\be^2 > \al^2$. 
Then $F_{\mu,\al,\be}$ takes values in the $4$-dimensional 
subgroup $(\C^{*} \times \GL(2,\C)) \cap \SL(3,\C)$,  
in particular, it is $\SL(2,\C)$-valued if $\mu=0$. 

For $F=F_{\mu,\al,\be}$, the $1$-form $\alpha_F$ 
is computed as follows:
\[
  \alpha_F =
    \begin{bmatrix} 
   (\mu +\frac{\al\be -3\mu^2}{\al + \be})z^{-1} & 0 & 
   -\frac{\sqrt{(\be^2+3\mu^2)(\al^2+3\mu^2)}}{\al+\be} z^{-\al-\be-1} \\
                            0   & -2\mu z^{-1} & 0 \\ 
  \frac{\sqrt{(\be^2+3\mu^2)(\al^2+3\mu^2)}}{\al+\be} z^{\al+\be-1} & 
  0 & (\mu -\frac{\al\be -3\mu^2}{\al + \be})z^{-1}
     \end{bmatrix}dz. 
\]
We treat the case when $f$
itself is single-valued on $\C \setminus \{ 0 \}$, i.e., 
\begin{equation*}
   f = F_{\mu, a,b}(F_{\mu,a,b})^*, \qquad 
        b-a \in \Z \setminus \{0 \},
\end{equation*}
where $F_{\mu,a,b}$ is given by \eqref{3:eq:Fmuab}. 
Then by \eqref{SL-killing}, it is verified that 
the induced metric $ds^2$ is non-degenerate and complete on $M$.
Moreover, $\alpha_F^\#=dFF^{-1}$ is computed as 
\begin{equation*}
  \alpha_F^\# = 
  \begin{bmatrix}
  \left( \mu + \frac{ab+3\mu^2}{b-a} \right)z^{-1} & 0 & -(a+b)z^{a-b-1} \\
  0 & -2 \mu z^{-1} & 0 \\
  (a+b)z^{b-a-1} & 0 & \left( \mu - \frac{ab+3\mu^2}{b-a} \right)z^{-1}
  \end{bmatrix}dz.
\end{equation*}
Indeed, $\alpha_F^\#$ is single-valued on $\C \setminus \{ 0 \}$. 
Furthermore, ${ds^2}^\#=\trace \alpha^\#(\alpha^\#)^*$ is 
\begin{equation*}
  {ds^2}^\# = \left\{
    \frac{2(a^2+3\mu^2)(b^2+3\mu^2)}{(b-a)^2}
    +(a+b)^2(|z|^{2a-2b}+|z|^{2b-2a})
     \right\}|z|^{-2} \,dz\, d\zb. 
\end{equation*}
  If $|a-b|=1$, the dual total curvature of $f$ is equal to $-4\pi$
  which satisfies equality in the Chern-Osserman type inequality
  (Theorem~\ref{main}) in the next section.
\end{example}
\section{Chern-Osserman type inequality}\label{sec:osserman}
Let $\Delta^{*}=\{z\in\C\,;\,0<|z|<1\}$ denote the unit disk punctured
at the origin. 
Let us consider a conformal immersion $f\colon{}\Delta^{*}\to N=G/H$ 
which has holomorphic right Gauss map, and
$F\colon\widetilde \Delta^* \to G$ its holomorphic lift, 
where $\widetilde \Delta^*$ is the universal cover of $\Delta^*$.
We denote by $f^\#$ the dual of $f$.
Recall that a metric on $\Delta^{*}$ is said to be {\it complete\/} at 
the origin if any path convergent to the origin has infinite length. 
If the metric $ds^2$ is asymptotic to the metric
\[
   |z|^{2\mu}dz\, d\bar z\qquad (\mu\in \R),
\]
then we call $\mu$ the order of the metric $ds^2$ at the
origin and denote by
\[
   \ord_{z=0}ds^2=\mu.
\]
It can be easily seen that the metric $ds^2$ is complete if and only if
$\ord_{z=0}ds^2\le -1$.

The following lemma is a generalization of the corresponding result
\cite{Yu} for \cmcone{} surfaces in $\Hyp^3$.
\begin{lemma}\label{lem:complete}
  The induced metric $ds^2$ of $f$ is complete at the origin if and 
  only if the induced metric ${ds^2}^\#$ of the dual $f^{\#}$ is also. 
\end{lemma}
\begin{proof}
Since $(ds^2{}^\#){}^\# = ds^2$, it suffices to show the one 
direction.
Let $\varGamma$ be a path tending to $0$ in $\Delta^{*}$. 
Under the assumption that the length $L^\#(\varGamma)$ 
of $\varGamma$ with respect to ${ds^2}^\#$ is finite, 
it is sufficient to prove $L(\varGamma) < \infty$. 

Denoting by $\widetilde{\varGamma}$ one of the lifts of $\varGamma$ to 
the universal cover $\widetilde \Delta^{*}$ of $\Delta^{*}$, 
 we can see from the completeness of $G/H$ that 
$f^\#(\widetilde{\varGamma})$ is bounded in $G/H$. 
The compactness of the fiber $H$ of the fiber bundle 
$G \to G/H$ implies that $F^\#(\widetilde{\varGamma})$ is bounded 
in $G$. 
Hence $F(\widetilde{\varGamma})$ is also bounded in $G$ because 
it is the image of $F^\#(\widetilde{\varGamma})$ under the 
diffeomorphism $a \mapsto a^{-1}$ of $G$. 
If we write $F=(F_{ij})$, then 
\begin{equation*}
  | F_{ij} | \le A \qquad\text{on}\qquad \widetilde{\varGamma}
\end{equation*}
holds for some constant $A$. 
Since 
\begin{equation*}
  \alpha_F=F^{-1}dF=F^{-1}(dFF^{-1})F=-F^{-1}\alpha_F^{\#}F,
\end{equation*}
we have 
\begin{equation}\label{estial}
  |\alpha_j|\leq C |\alpha_j^{\#}|\qquad
       (j=1,2,\dots,N)
\end{equation}
for some constant $C$. 
Thus, we have
\begin{equation*}
 L(\varGamma)=\int_{\varGamma} ds
             =\int_{\varGamma}\sqrt{\sum_{j=1}^n|\alpha_j|^2}
             \leq C\int_{\varGamma}\sqrt{\sum_{j=1}^n|\alpha_j^{\#}|^2}
             =CL^{\#}(\varGamma)<\infty \, .
\end{equation*}
\end{proof}
Let $f\colon{}\Delta^*\to N$ be a complete conformal immersion with
holomorphic right Gauss map.
Then, by Lemma~\ref{lem:complete}, the metric ${ds^2}^\#$ is complete at the 
origin if $f\colon{}\Delta^*\to N$ is a complete immersion.
Note that the assertion is intrinsic for $(M,{ds^2}^\#)$. 
Hence, in the similar way to that of minimal surface theory,
we have the following Lemma (see \cite{lawson}).
\begin{lemma}\label{lem:os2}
  Suppose $f$ is complete at the origin and 
  has finite dual total curvature i.e., 
\begin{equation}
   \int_{\Delta^{*}}(-K^{\#})\,dA^{\#} < \infty.
\end{equation}
  Then $\alpha_F^{\#}$ has a pole at the origin. 
\end{lemma}
Now we shall state our main result,
which is a generalization of the case of \cmcone{} surface 
in $\Hyp^3$ (\cite{UY4} and \cite{Yu2}) 
and the original Chern-Osserman type inequality for
minimal surfaces in $\mathbf R^r(\subset G)$
where $r$ is the dimension of maximal abelian 
subgroup of $G$. 

\begin{theorem}[Chern-Osserman type inequality]\label{main}
  Let $f\colon{}M\to G/H$ be a complete conformal  
  immersion with holomorphic right Gauss map
  of finite dual total curvature, and let  $f^{\#}$ be its dual. 
  Then the total dual curvature satisfies 
\[
   \frac{1}{2\pi}\int_M (-K^{\#})\,dA^{\#}\geq
        -\chi(M)+(\text{\rm{the number of ends}}). 
\]
  Here $\chi(M)$ denotes the Euler number of $M$. 
\end{theorem}
The following lemma is crucial to prove the theorem: 
\begin{lemma}\label{lem:os3}
  Let  $f\colon{}\Delta^{*}\to G/H$ be a conformal  
  immersion which has holomorphic right Gauss map. 
  Then the order of the metric ${ds^2}^\#$ of $f^\#$
  at the origin is less than or equal to  $-2$.
\end{lemma}
First, we shall prove Theorem~\ref{main} using Lemma~\ref{lem:os3}.
\begin{proof}[Proof of Theorem~{\rm\ref{main}}]
  Let $F$ be a holomorphic lift of $F$ and $\alpha^{\#}=dFF^{-1}$.
  Using the basis \eqref{eq:basis} of $\Lie{g}$, 
  we write $\alpha^{\#}=\sum_{j=1}^n \alpha_j^{\#} e_j$, and define a $\C^n$-valued
  $1$-form
\[
    \alpha^{\#}_0 = \left(\alpha^{\#}_1,\dots,\alpha^{\#}_n\right)
\]
  on $M$.
  Let $f^{\#}_0 \colon \widetilde M \to \R^n$ be a 
  conformal minimal immersion defined by 
\[
  f^{\#}_0 = 2 \Re \int_{z_0}^z \alpha^{\#}_0.
\]
  Then the Gauss map of $f^\#_0$ is, by definition, 
\[
   \nu_0^\#= [ \alpha_1^\#: \cdots : \alpha_n^\#] \colon
         \overline M \to {\mathcal Q}_{n-2} \subset \CP^{n-1},
\]
  where $\overline M$ is the compactification of $M$
  (see Proposition \ref{prop:cp-fi}).
  Since $\alpha^{\#}$ is single-valued on $M$ by
  Proposition~\ref{prop:metric}, $\nu_0^\#$ is a single-valued
  map defined on $\overline M$. 
  Since the induced metric $ds^2_0{}^{\#}$ of $f^{\#}_0$ is given 
  by 
\[
    ds^2_0{}^{\#} = \sum_{j=1}^n \alpha_k^{\#}\overline\alpha_k^{\#}
\]
  (see \eqref{eq:min-met}), 
  we have ${ds^2}^{\#}={ds^2_0}^{\#}$.
  Hence the dual total curvature
\[
    \int_M (-K^{\#})\,dA^{\#}
\]
  is the total curvature of the minimal immersion $f_0^\#$.
  As in the proof of Lemma~11 of \cite{lawson}, we have
\begin{equation}\label{K1}
  \int_M (-K^{\#})\,dA^{\#} = 2\pi k,
\end{equation}
  where $k$ is the homology degree of $\nu_0^\#$.
  Let $p_1,\dots,p_r\in \overline M$ be the ends of $f_0^\#$. 
  By Lemma~\ref{lem:os2}, $\alpha_j^{\#}$ has a pole at each end.
  Let $m_j$ denote the maximum order of the pole of 
  $\alpha_1^{\#},\dots,\alpha_n^{\#}$ at the end $p_j$.
  Then, as in the proof of Theorem~19 of \cite{lawson}, we have
\begin{equation}\label{K2}
    k=\sum_{j=1}^r m_j-\chi(\overline{M}),
\end{equation}
  where $r$ is the number of ends. 
  It should be remarked that to prove \eqref{K1}
  and \eqref{K2}, we do not need well-definedness of $f_0^{\#}$ on $M$, but
  only the fact that the Gauss map $\nu_0^{\#}$
  is single-valued on $M$.
  It can be easily checked that the order of ${ds^2}^\#$
  at $p_j$ is equal to $-m_j$.
  By Lemma \ref{lem:os3}, we have
\begin{equation}\label{K3}
     m_j\ge 2.
\end{equation}
  By  \eqref{K1}--\eqref{K3}, we have
\[
   \frac 1{2\pi}\int_M (-K^{\#})\,dA^{\#} =k = \sum_{j=1}^r m_j-\chi(\overline M) 
                 \ge 2r-2\chi(\overline{M})=r-\chi(M).
\]
  This proves the theorem.
\end{proof}
\begin{remark}\label{rem:finn}
  There is an alternative proof of the theorem: 
  Finn \cite{f} (Theorem~20) proved the equality
\begin{equation}\label{eq:finn}
    \frac1{2\pi}\int_{M} K \,dA=\chi(M)-\sum_{j=1}^r t_j
\end{equation}
  for a certain class of complete Riemannian $2$-manifolds $M$
  that are biholomorphic to punctured closed Riemann surfaces
  $\overline M\setminus \{p_1,\dots,p_r\}$.
  For each punctured point $p_j$, we can take a complex coordinate
  $z$ ($|z|>R$) which maps $\{|z|>R\}$ to the punctured neighborhood
  of $p_j$ in $M$. Then $t_j$ is defined by
\[
   t_j=\lim_{r\to \infty} \frac{\mathcal L(r)^2}{4\pi \mathcal A(r;R)},
\]
  where $\mathcal L(r)$ is the length of the curve $|z|=r$ and
  $\mathcal A(r;R)$ is the area of the region $R<|z|<r$.
  Applying this formula for $(M,{ds^2}^\#)$.
  It can be easily checked that
\[
    t_j=-\ord_{z=p_j}{ds^2}^\#.
\]
  Thus Theorem \ref{main} follows from the fact $t_j\geq 2$,
  shown in Lemma~\ref{lem:os3}.
  
\end{remark}

The rest of this section is devoted to proving Lemma \ref{lem:os3}. 
Let  $f\colon{}\Delta^{*}\to G/H$ be  a conformal immersion which has
holomorphic right Gauss map and $F\colon{}\widetilde{\Delta^{*}}\to G$ 
a holomorphic lift of $f$. 
We consider $F$ to have values in $\SL(n,\C)$.
We set $\alpha^\#=-dFF^{-1}$.

Then the column vectors of $F$ consist of a fundamental system of
solutions of the ordinary differential equation
$dy=-\alpha^\# y$, where $y$ is a $\C^n$-valued function. 
Then it can be easily checked that $FF^*$ is single-valued on $\Delta^*$
if and only if the monodromy group  of the differential equation 
$dF=-\alpha^\# F$ is contained in $H$.
Moreover, $-\ord_{z=0}{ds^2}^\#$ is equal to
the order of pole of $\alpha^\#$ at the origin.
Recall that $G$ is a connected and semi-simple subgroup of $\SL(n,\C)$
($n=\dim G$) and $H\subset \SU(n)$.
Thus by setting $\beta=-\alpha^\#$, 
Lemma~\ref{lem:os2} is a corollary of the following: 
\begin{assertion}\label{ass2}
  Let $\beta$ be an ${\mathfrak {sl}}(n,\C)$-valued null 
  holomorphic $1$-form on $\Delta^{*}$. Suppose that it has a 
  pole of order $1$ at the origin. Then the monodromy group 
  of the differential equation $dFF^{-1}=\beta$
  is not contained in $\SU(n)$. 
\end{assertion}
\begin{proof}
  Suppose that the monodromy group is included in $\SU(n)$. 
  We shall make a contradiction.
\paragraph{Step 1.}
  First, we shall describe the monodromy matrix. 
  Write $\beta=b(z)dz$ and consider a system of ordinary 
  differential equations 
\begin{equation}\label{eq:ode0}
    \frac{dy}{dz}=b(z) y 
\end{equation}
  for $y=y(z)\in \C^{n}$. 
  Take a base point $z_0\in \Delta^{*}$.
  Let $F_0\colon{}\Delta^*\to \SL(n,\C)$ ($F_0(z_0)=\id$) 
  such that the column vectors of $F_0$ are a fundamental system 
  of solutions of \eqref{eq:ode0}.  

  Let us denote the monodromy representation by $\rho_0$, i.e., 
  $\rho_0$ is determined by 
  $ F_0\circ \tau = F_0\rho_0(\tau)$ 
  for any covering transformation $\tau$ of $\widetilde \Delta^*$. 
  By our assumption, $\rho_0(\tau)\in \SU(n)$, in particular, 
  eigenvalues of $\rho_0(\tau)$ are of the 
  form $e^{i\theta}$ ($\theta\in\R$).

  Since $F_0$ is non-singular at each point of $\Delta^*$,
  an arbitrary fundamental system $F$ of 
  solutions to \eqref{eq:ode0} can be written as 
\[
   F=F_0 P, \qquad P\in\GL(n,\C).
\]
  We let $\rho(\tau)$ the monodromy matrix
  determined by $ F\circ  \tau= F \rho(\tau)$. 
  Then we can write
\begin{equation}\label{eq:mono}
  F\circ \tau = F \rho(\tau),\qquad \rho(\tau)=P^{-1}\rho_0(\tau)P. 
\end{equation}
  Hence, eigenvalues of $Q$ are 
  also of the form $e^{i\theta}$ ($\theta\in\R$) and $\rho(\tau)$ is
  diagonalizable.
\subparagraph{Step 2}
  Since we assumed that $\beta=b(z)dz$ has a pole of order $1$ at the
  origin, $b(z)$ has the Laurent expansion 
\[
    b(z)=\frac{1}{z}R+ \sum_{j=0}^{\infty} z^j A_j
    \qquad (R\neq 0).
\]
 We shall show that all eigenvalues of $R$ vanish.
 Let $\lambda$ be an eigenvalue of $R$. 
 First we shall prove that $\lambda$ is a real number.
 We have only to prove it for an eigenvalue $\lambda$ 
 such that $\lambda + j$ is not an eigenvalue for any 
 positive integer $j$. 
 (In fact, whenever both $\lambda$ and $\lambda +j$ are eigenvalues for
 an integer $j$, $\lambda$ is real if and only if $\lambda +j$ is also.)
 By Theorem~\ref{ap:th:3} in Appendix, 
 there exists a solution $y$ to \eqref{eq:ode0} 
 which is expanded as
\[
    y=z^{\lambda}\left[ v_0 + z v_1+z^2 v_2 + \cdots\right] \qquad 
       (v_0\neq 0).
\]
  Denote by $\tau$ the covering transformation corresponding to 
  a loop going once around the origin in $\Delta^{*}$ in 
  the positive direction. 
  Then $y$ satisfies 
\[
    y\circ \tau = e^{2\pi i\lambda} y.
\]
  If we take a fundamental system $F=[y_1, y_2,\dots, y_n]$ 
  of solutions to \eqref{eq:ode0} so that the first column vector 
  $y_1$ is equal to $y$, then 
\[
   F\circ \tau =
      [e^{2\pi i\lambda}y,y_2\circ \tau, \dots,
                                        y_n\circ \tau]
    =[y, y_2,\dots, y_n]
      \left[
       \begin{array}{cccc}
            e^{2\pi i\lambda} & * & \dots & * \\
             0 & * & \dots & * \\
             \vdots &  & \ddots &  \\
             0 & * & \dots & *
        \end{array}
      \right]. 
\]
  Hence, $e^{2\pi i \lambda}$ is an eigenvalue of the 
  monodromy matrix $\rho(\tau)$ of \eqref{eq:mono}. 
  Then it follows from Step 1 that $\lambda$ must be real. 
  Thus we can conclude  that all the eigenvalues of $R$ are real. 
  Moreover, 
  the nullity condition $B(\beta,\beta)=0$ gives 
\[
    0=\trace b(z)^2= \frac{1}{z^2}\trace R^2+\cdots, 
\]
  in particular, $\trace R^2=0$. 
  Namely, the square-sum of eigenvalues of $R$ is zero.
  Hence, all the eigenvalues of $R$ must be zero. 

\paragraph{Step 3}
  Finally, we shall make a contradiction.
  Suppose  $R$ is diagonalizable.
  Then $R$ vanishes, which contradicts the completeness around the
  origin. 
  Hence, $R$ is not diagonalizable, and its Jordan normal form 
  is 
\[
  \left[
   \begin{array}{ccccc}
   0 & 1 &   & &  \\
     & 0 & * & \smash{\lower-0.5ex\hbox{\Huge 0}}&       \\
     &   & \ddots & \ddots & \\
     & \smash{\lower2ex\hbox{\Huge 0}} & & \ddots & * \\
     &                         & &        & 0
  \end{array}\right], 
\]
  where $*$ denotes $0$ or $1$. 
  By Theorem~\ref{ap:th:4} of Appendix, there exist linearly independent 
  solutions $y$, $\tilde y$ of \eqref{eq:ode0} such that 
\begin{align*}
    y&= p_0+zp_1+z^2 p_2+\cdots\\
    \tilde y&=(q_0+zq_1+z^2 q_2+\cdots)+(\log z) y.
\end{align*}
  For the covering transformation $\tau$ as above, 
\[
    y\circ \tau=y,\qquad
    \tilde y\circ \tau= \tilde y+2\pi i y. 
\]
  If we take a fundamental system $F$ of solutions so that 
\[
     F=[y,\tilde y,y_3,\dots,y_n], 
\]
  the monodromy matrix is given by
\[
  F\circ \tau = F
              \left[
               \begin{array}{ccccc}
                1 & 2\pi i & * & \dots & * \\
                0 &   1    & * & \dots & * \\
                0 &   0    &   &       &   \\
                \vdots & \vdots & & \smash{\lower2ex\hbox{\Huge *}} & \\
                0 &   0    &   &       &    
               \end{array}
              \right].
\]
  This means that the monodromy matrix $\rho(\tau)$ of \eqref{eq:mono} is
  not diagonalizable, which contradicts the conclusion of Step 1.
\end{proof}

\begin{remark}
 Consider a complete conformal immersion $f$ of 
 $M=\overline{M}\setminus\{p_1,\dots,p_r\}$ into $G/H$
 with holomorphic right Gauss map.
 Then completeness of the induced metric $ds^2$ 
 at each end $p_j$ of $f$  implies that
\begin{equation}\label{eq:non-dual-osserman}
    u_j := -\ord_{z=p_j}{ds^2}\geq 1.
\end{equation}
By \eqref{eq:finn} in Remark~\ref{rem:finn},
 the total absolute curvature of $f$  satisfies
\[
   \frac{1}{2\pi}\int_M (-K)\,dA \geq -\chi(\overline M)+ r 
                                 =\chi(M),
\]
 which is the well known Cohn-Vossen inequality.

 When $N=\SL(2,\C)/\SU(2)=\Hyp^3$, we have $u_j>1$, 
 and consequently the equality in the Cohn-Vossen
 inequality does not hold (\cite{UY1}).
  For a general ambient space $G/H$, it is still unknown whether 
 $u_j>1$ holds or not.
 The technique as in the proof of Assertion~\ref{ass2}  cannot
 be applied for this problem, because the coefficient 
 matrix $\beta(z)$  might not be single-valued on $\Delta^*$.
\end{remark}

\section{Perturbation of minimal surfaces}\label{sec:perturb}
In Section~\ref{sec:gauss}, we have a non-trivial
example of surfaces with holomorphic right Gauss map in
$\SL(3,\C)/\SU(3)$. 
Perturbing minimal surfaces in $\R^n$,
we shall give in this section a further example of complete
surface in $\SL(3,\C)/\SU(3)$ with holomorphic right Gauss map,
which satisfy the equality of the Chern-Osserman type inequality
(Theorem~\ref{main}).
The method is similar to that in \cite{UY2,ruy1}, in which
many examples of \cmcone{} surfaces in $\Hyp^3$ are obtained
as a perturbation of minimal surfaces in $\R^3$.
Thus the method in this section will provide many examples
of surfaces with holomorphic right Gauss map.

\subsection*{Perturbation of monodromy matrices}
Let $M$ be a Riemann surface and assume that 
the fundamental group $\pi_1(M)$ of $M$ 
is finitely generated by $\{\gamma_1,\dots,\gamma_k\}$, 
where $\gamma_j$ ($j=1,\ldots,k$) is a loop based at 
$z_0 \in M$. 
Let $\hat\pi \colon \widetilde M \to M$ be the universal 
covering of $M$ and take $\tilde z_0 \in \hat\pi^{-1}(z_0)$. 
Denote by $\tilde\gamma_j$ the lift of $\gamma_j$ emanating 
from $\tilde z_0$, that is, 
\[
  \tilde\gamma_j 
     = \gamma_j \circ \hat\pi\colon{}[0,1]
     \longrightarrow \widetilde M,
\qquad
     \tilde\gamma_j(z_0)=\tilde z_0.
\]

A  $\C^n$-valued holomorphic $1$-form 
$\alpha_0=(\alpha_1,\dots,\alpha_n)$ on $M$ is called the Weierstrass data
if it satisfies 
\[
  \sum_{j=1}^n\alpha_j\cdot\alpha_j=0\qquad\text{and}\qquad
  \sum_{j=1}^n\alpha_j\cdot\overline{\alpha_j}>0.
\]
A Weierstrass data $\alpha_0$ induces a conformal minimal immersion
\[
   f=2\Re\int_{z_0}^{z}\alpha_0 \colon \widetilde M\longrightarrow \R^n.
\]
Here $f$ is single-valued on $M$ if and only if
\[
   \Re\int_{\gamma_j}\alpha_0=0\qquad (j=1,\dots,k).
\]
Let $G\subset \SL(n,\C)$ be a Lie group of dimension $n$ as in
Section~\ref{sec:gauss} and set
\[
   \hat\alpha := \sum_{j=1}^n \alpha_j e_j,
\]
where $\{e_1,\dots,e_n\}$ is a basis of $\Lie{g}$ as in 
\eqref{eq:basis}.
Note that $\alpha_0$ is null if and only if $\hat\alpha$ is also. 
Let us consider the initial value problem of the differential 
equation 
\begin{equation}\label{eq:5:ode}
    dF_{c}=-c\hat\alpha F_{c},\qquad
    F_{c}(\tilde z_0)=\id 
\end{equation}
on $\widetilde M$, where $c$ is a real constant. 
For the (unique) solution $F_{c}\colon \widetilde M\to G$ 
of \eqref{eq:5:ode}, we define a conformal immersion
with holomorphic Gauss map
\[
    f_{c}=F_{c}(F_{c})^*. 
\]
For $j=1,\dots,k$, we let $\rho^j({c})$ be the monodromy
matrix defined by
\[
  F_{c}\circ \tau_j=F_{c} \rho^j(c),
\]
where $\tau_j$ is the covering transformation
such that $\tau_j(\tilde z_0)=\tilde\gamma(1)$. 
We set
\[
  \sigma^j(c)=\rho^j({c})(\rho^j({c}))^*. 
\]
Note that $f_{c}$ is well-defined on $M$ if and only if 
$\sigma^j(c)$ is identity for every $j=1,\dots,k$. 

\begin{lemma}\label{6:th:perturb}
  Under the notations above, 
  we assume  $f_0\colon{}M\to\R^n$ is single-valued on $M$.
  Then the following identity holds{\rm :}
\[
  \sigma^j(c) = \id -2 c\sum_{k=1}^n \left(\Re \int_{\gamma_j} 
   \alpha_k\right)
   e_k + o(c). 
\] 
\end{lemma}
\begin{proof}
  Since $G$ is a subgroup of $\SL(n,\C)$, the Cartan model of 
  $G/H$ is included in the set of Hermitian matrices. 
  So each $e_j$ is Hermitian:
\[
     e_{j}^*=e_j.
\]
  Since $F_{0}$ is identity, we have
\begin{equation}\label{6:condi:wd}
  \rho^j({0})=\id \qquad (j=1,2,\dots, k).
\end{equation}

  Differentiating both sides of \eqref{eq:5:ode} with respect to 
  $c$ at $c=0$, we have 
\[
   d\left.\frac{\partial F_{c}}{\partial c}\right|_{c=0}
      =-\hat\alpha F_{0}
      =-\hat\alpha \qquad (F_{0}=\id). 
\]
It follows that 
\[ 
    \left.\rho^j({c})\right|_{c=0}=
    \left.\frac{\partial F_{c}}{\partial c}(\tilde\gamma_j(1))
     \right|_{c= 0}
     = \oint_{\gamma_j} (-\hat\alpha)
     = -\sum_{i=1}^{n}\left[\oint_{\gamma_j}
     \alpha^i\right]e_i. 
\]
Using the fact that $e_k$ is Hermitian, we obtain 
\begin{align*}
  \left.\frac{\partial}{\partial c}\right|_{c=0}
   \sigma^j(c)&=
  \left.\rho^j({c})\right|_{c=0}\cdot
   {}^t\overline{\rho^j(0)}+
   {\rho^j({0})}\cdot
   \left.{}^t\overline{\rho^j({c})}\right|_{c=0}\\
  &=
  \left.\rho^j({c})\right|_{c=0}+
   \left.{}^t\overline{\rho^j({c})}\right|_{c=0}
   = -\sum_{i=1}^{n}\left[\oint_{\gamma_j}\alpha^i\right]e_i
     -\sum_{i=1}^{n}\left[\oint_{\gamma_j}\overline{\alpha}^i
         \right]{}^{t}\overline{e_i}&\\
   &= -\sum_{i=1}^{n} \left(\Re\oint_{\gamma_j}2\alpha^i\right) e_i
 \qquad &\qedsymbol
\end{align*}
\renewcommand{\qedsymbol}{\relax}
\end{proof}
\subsection*{An example}
We construct an example of genus zero with three ends
in the symmetric space $\SL(3,\C)/\SU(3)$. 
Let us consider a minimal surface 
\[
   f_0\colon{}\C\setminus\{0,-1\}\longrightarrow \R^4
\]
defined by 
\begin{equation}\label{eq:weier}
   f_0=2\Re\int_{z_0}^z \alpha_0,\qquad
   \alpha_0=
   (1- \langle g_0,g_0\rangle,i(1+\langle g_0, g_0\rangle), 2g_0)\omega_0,
\end{equation}
where $(g_0,\omega_0)$ is 
\begin{equation}\label{eq:w-data}
   g_0=\left(\frac{z+1}{z-1},2\frac{(z+1)^2}{(z-1)^2}\right),\qquad
   \omega_0 = \frac{(z-1)^4}{z^2(z+1)^2}\,dz. 
\end{equation}
It is easily verified $f_0$ is a single valued
minimal immersion on $\C\setminus\{0,-1\}$. 
It can be also checked that the induced metric 
\[
   ds^2_0=(1+2\langle g_0,\bar g_0\rangle+
            \langle g_0, g_0\rangle 
       \langle\bar g_0,\bar g_0\rangle)\omega_0\bar\omega_0
\]
is positive definite on $\C\setminus\{0,-1\}$ and complete.

Let us consider a seven-parameter family of $(G,\omega)$ 
defined by 
\[
   \alpha=
   (1- \langle g, g\rangle,i(1+ \langle g, g\rangle), 2g)\omega
\]
\begin{equation}\label{eq:w-deform}
 g=\left( \frac{a_2 z + a_3}{z-a_1},
          \frac{a_4 z^2 + a_5z+a_6}{(z-a_1)^2}\right), \qquad
  \omega=a_7\frac{(z-a_1)^4}{z^2(z+1)^2}\,dz.
\end{equation}
Here, seven parameters $\bmath{a}=(a_1,a_2,a_3,a_4,a_5,a_6,a_7)$ 
can run over $\C$. In other words, \eqref{eq:w-deform} is 
parameterized by fourteen real numbers. 
Note that $(g,\omega)$ at $\bmath{a}_0=(1,1,1,2,4,2,1)$ is the 
data \eqref{eq:w-data} of the initial surface $f_0$. 
It can be checked that 
\[
   ds^2=(1+2\langle g,\bar g\rangle+
           \langle g, g\rangle \langle\bar g,\bar g\rangle)
           \omega\bar\omega
\]
is a complete metric if $\bmath{a}$ is sufficiently close to
$\bmath{a}_0$. 

Now we set
\[
   G=\SL(3,\C),\qquad H=\SU(3)\qquad\text{and}\qquad N=G/H .
\]
Then 
\begin{align*}
   &e_1=\frac{1}{2\sqrt{3}}\left[\begin{array}{ccc}
             0 & 0 & 1 \\ 0& 0 & 0 \\ 1 & 0 & 0 
             \end{array}\right],\qquad
   e_2=\frac{1}{2\sqrt{3}}\left[\begin{array}{ccr}
             0 & 0 & -i \\ 0& 0 & 0 \\ i & 0 & 0 
             \end{array}\right], \\
   &e_3=\frac{1}{2\sqrt{3}}\left[\begin{array}{crc}
             1 & 0 & 0 \\ 0& -1 & 0 \\ 0 & 0 & 0 
             \end{array}\right],\qquad
   e_4=\frac{1}{6}\left[\begin{array}{ccr}
             1 & 0 & 0 \\ 0& 1 & 0 \\ 0 & 0 & -2 
             \end{array}\right],\qquad
\end{align*}
are orthonormal in $T_oN$.
We set 
\[
  \hat\alpha=\alpha_1 e_1 + \alpha_2 e_2+
         \alpha_3 e_3 + \alpha_4 e_4.
\]
Let $G'$ be a subgroup of $\SL(3,\C)$ defined by
\[
    G'=\{a\in \SL(3,\C)\,;\, a_{12}=a_{21}=a_{23}=a_{32}=0\}, 
\]
where $a_{ij}$ is the $(i,j)$-component of the matrix $a$. 
Then $\hat\alpha$ is valued in the Lie algebra of $G'$.

Let $\gamma_1$ and  $\gamma_2$ be loops on $\C\setminus\{0,1\}$
surrounding $0$ and $1$ respectively.
Then $\{\gamma_1,\gamma_2\}$ generates the fundamental group of
$\C\setminus\{0,1\}$.
We use the same notation in the first part of this section.
The solution $F_c$ of \eqref{eq:5:ode} lies in $G'$, in particular
we have
\begin{equation}\label{eq:vanish}
  \sigma^j_{12}(c)=\sigma^j_{21}(c)
        =\sigma^j_{23}(c)=\sigma^j_{32}(c)=0\qquad (j=1,2),
\end{equation}
where $\rho^j(c)=(\rho^j_{kl}(c))_{k,l=1,2,3}$.
We set
\begin{multline*}
  \varphi(c,\bmath{a})=\left(c,
  \frac{\sigma^1_{11}(c)-1}c,
    \frac{\sigma^1_{33}(c)-1}c,
    \frac{\sigma^1_{13}(c)}c,\frac{\sigma^1_{31}(c)}c,\right. \\
  \left.\qquad \frac{\sigma^2_{11}(c)-1}c,
    \frac{\sigma^2_{33}(c)-1}c,
    \frac{\sigma^2_{13}(c)}c,\frac{\sigma^2_{31}(c)}c\right).
\end{multline*}
Then by Lemma~\ref{6:th:perturb}, $\varphi$ is a smooth map from 
a neighborhood of $(c,\bmath{a}_0)\in\R^{15}$ into 
$\R^9$ which is written as
\begin{multline*}
   \varphi(c,\bmath{a})= 
    (c,\bmath{0})\\
    -4\pi\left(0,
         \Im\Res_{z=0}\alpha_1,\dots,
         \Im\Res_{z=0}\alpha_4,
         \Im\Res_{z=-1}\alpha_1,\dots,
         \Im\Res_{z=-1}\alpha_4,\right)+o(c).
\end{multline*}
Then the Jacobian matrix $J$ of $\varphi$ at
$(c,\bmath{a})=(0,\bmath{a}_0)$ is in the form
\[
    J = \begin{bmatrix} 
          1 & *\\
          \bmath{0} & 4\pi{
          \dfrac{\partial\left(\Im\dRes_{z=0} \alpha, \Im\dRes_{z=-1}\alpha\right)}{
          \partial (\Re a_1,\Im a_1,\dots,\Re a_7,\Im a_7)}}
        \end{bmatrix}.
\]
Here, by direct calculation, we have 
\begin{equation}\label{6:eq:jac}
  \rank \left.\frac{\partial(
         \Im\Res_{z=0}\alpha, \Im\Res_{z=1}\alpha)}{
    \partial(\Re a_1,\Im a_1,
             \Re a_2,\Im a_2,\dots,
             \Re a_7,\Im a_7)} \right|_{\smath{a}_0}=8. 
\end{equation}
Hence there exists $\bmath{a}=\bmath{a}(c)$ such that
$\varphi(c,\bmath{a}(c))=0$. 
Then this implies
\[
  \sigma^1(c)=\sigma^2(c)=\id.
\]
Since $\rho^3(c)=(\rho^1(c)\rho^2(c))^{-1}$, we have
$\sigma^3(c)=\id$. Therefore we have conformal immersions 
with holomorphic right Gauss map
$f_c \colon \C\setminus\{0,-1\} \to G/H$.

Finally, we calculate the total dual absolute curvature of 
our example $f_c$. 
The total absolute curvature of the initial surface $f_0$ is 
an integral multiple of $2\pi$ and the total dual curvature 
of $f_c$ is smoothly depending on the parameter $c$.
Hence they coincide.
Since the total absolute curvature of $f_0$ is $8\pi$, 
we have the total absolute curvature of $f_c^{\#}$ is also $8\pi$,
which satisfies equality in the Chern-Osserman type inequality of 
Theorem~\ref{main}.  
\appendix
\begingroup
\renewcommand{\thesection}{\Alph{section}}
\section{}
\label{app:A}
\endgroup
\small
We review some facts from the general theory of ordinary 
differential equations from \cite{CL}
which are needed in the proof of the 
local Chern-Osserman type inequality (see Section~\ref{sec:osserman}).
Let $A(z)$ be a $\Lie{gl}(n,\C)$-valued holomorphic 
function on $\Delta^*=\{z\in\C\,;\, 0<|z|<1\}$. 
Assume that $A(z)$ has a pole at $z=0$. We consider the following 
linear ordinary equation: 
\begin{equation}\label{eq:ode}
   \frac{dy}{dz}=A(z)y
\end{equation}
for an unknown $\C^n$-valued (column vector-valued) function $y$. 
The fundamental system of solutions to \eqref{eq:ode} is, by definition, 
a matrix whose column vectors consist of 
$n$ linearly independent solutions. 

\begin{fact}[\cite{CL}, Section 4]\label{fact:1}
  The fundamental system $\Phi(z)$ of solutions to \eqref{eq:ode} can be 
  written as 
\begin{equation}\label{eq:fundsys}
   \Phi(z)=S(z)z^{\Lambda}\qquad z^{\Lambda}=e^{\log z \cdot {\Lambda}}
\end{equation}
  where $S(z)$ is a holomorphic function on $\Delta^*$ and  
  ${\Lambda}$ is a constant matrix. 
\end{fact}
If we replace the fundamental system $\Phi(z)$ into 
$\widetilde \Phi (z)$ so that $\widetilde\Phi(z)=\Phi(z)P$ 
for some $P\in\GL(n,\C)$, then we have 
\[
   \widetilde\Phi(z)=S(z)z^{\Lambda}P=\{S(z)P\}P^{-1}z^{\Lambda}P
                             =\{S(z)P\}z^{P^{-1}{\Lambda}P}. 
\]
Hence, we may assume that the matrix $M$ in \eqref{eq:fundsys} 
is the Jordan normal form, without loss of generality.  

With respect to the equation \eqref{eq:ode}, 
the origin $0$ is called a {\it regular singularity\/} 
if $S(z)$ in the system of solutions 
\eqref{eq:fundsys} has at most a pole. 
If $S(z)$ has an essential singularity at $0$, 
it is called an {\it irregular singularity}. 

\begin{fact}[{\cite[Section 4]{CL}}]\label{fact:2}
 The origin $0$ is a regular singularity of the equation \eqref{eq:ode} 
 if $A(z)$ has a pole of order $1$ at the origin. 
\end{fact}

Hereafter we assume that $A(z)$ has a pole of order $1$ at the origin and  
the Laurent expansion of $A(z)$ is given by 
\[
  A(z)=\frac{1}{z}R+\sum_{j=0}^{\infty}z^jA_j
    \qquad\text{($R$, $A_j$ are constant matrices).}
\]
We denote by $\lambda_0$ the eigenvalue of $R$ such that 
$\lambda_0 + j$ is not an eigenvalue of $R$ for any positive integer $j$. 

\begin{theorem}\label{ap:th:3}
  Under the assumption above,  
  there exists a solution $y$ to the equation \eqref{eq:ode} 
  which is of the form 
\[
     y=z^{\lambda_0}\left[p_0+zp_1+z^2 p_2+ \cdots\right]. 
\]
\end{theorem}
\begin{proof}
Let us take an arbitrary constant vector $y_0\in\C^n$. 
For a complex parameter $\lambda$, we define 
$y_0(\lambda,y_0):=y_0$ and 
$y_1(\lambda,y_0)$, $y_2(\lambda,y_0)$, \dots 
inductively by 
\begin{equation}\label{eq:formal}
     y_{j+1}(\lambda,y_0)
     =-[(\lambda+j+1)I-R]^{-1}\sum_{k=0}^{j}A_ky_{j-k}(\lambda,y_0)
   \qquad(j=0,1,2,\dots)
\end{equation}
where $I$ is the identity matrix. 
The right-hand side of \eqref{eq:formal} implies that the
$y_j$'s are $\C^n$-valued rational functions in the variable $\lambda$, 
with a regular point $\lambda=\lambda_0$.  

We consider a formal power series 
\begin{equation}\label{eq:power}
    y(\lambda,y_0):=z^{\lambda}\sum_{j=0}^{\infty}z^j y_j(\lambda,y_0). 
\end{equation}
The ordinary differential equation \eqref{eq:ode} has
$n$ linearly independent solutions in the space of formal
power series of $\log z$, and by Facts~\ref{fact:1} and \ref{fact:2},
we know that these formal power series converge. 
So we may treat the solution as a formal power series.
Differentiating \eqref{eq:power}, we have 
\begin{multline}\label{eq:formal-diff} 
   \left[\frac{d}{dz}-A\right]y(\lambda,y_0)\\
     = z^{\lambda}\left[
            \frac{1}{z}(\lambda I-R)y_0+
            \sum_{j=0}^{\infty}
              z^j\left[
   \left\{(\lambda+j+1)I-R\right\}
                      y_{j+1} + \sum_{k=0}^{j}A_ky_{j-k}\right]
        \right]\\
    =z^{\lambda-1}(\lambda I-R) y_0
\end{multline}
By our assumption, $(\lambda_0+j+1)I-R$ is a nonsingular matrix. 
Hence $y_j(\lambda_0,y_0)$ defined in \eqref{eq:formal} 
are regular values (i.e., $\ne \infty$).  Moreover, if  
$p_0$ is a $\lambda_0$-eigenvector of $R$, then 
$y(\lambda_0,p_0)$ is a solution to the equation \eqref{eq:ode},
because $(\lambda_0 I-R)p_0=0$. 
\end{proof}
Furthermore, we assume that $\lambda_0$ is an eigenvalue 
of multiplicity $m \ge 2$. If the $\lambda_0$-eigenspace has 
the maximal dimension $m$, then there exist $m$ linearly 
independent solutions. Indeed, we can construct them by 
Theorem \ref{ap:th:3}. 
Let us consider the other case, that is, 
suppose that the Jordan canonical form of $R$ is given by 
\[
    P^{-1}RP=
       \left[\begin{array}{ccccc}
     \lambda_0 & 1 &   & &  \\
     & \lambda_0 & *&\smash{\lower-0.5ex\hbox{\Huge 0}}&       \\
     &   & \ddots & \ddots & \\
     & \smash{\lower2ex\hbox{\Huge 0}} & & \ddots & * \\
     &                         & &        & *
    \end{array}\right],\qquad
   (P\in\SL(n,\C)). 
\]
\begin{theorem}\label{ap:th:4} 
  Under the assumption above, 
  there exist two linearly independent 
  solutions of the equation \eqref{eq:ode}, which are 
  of the form  
\begin{align*}
     y&=z^{\lambda_0}\left[p_0+zp_1+p^2 y_2+ \cdots\right]\\
     \tilde y&=z^{\lambda_0}\left[q_0+zq_1+z^2 q_2+ \cdots\right]+
             (\log z)y. 
\end{align*}
\end{theorem}
\begin{proof}
  Let $p_0$ be the first column of $P$ and $q_0$ the second column.
  Then we have
\begin{equation}\label{eq:eigen-vector}
     Rp_0=\lambda_0p_0\qquad
     Rq_0=\lambda_0q_0+p_0. 
\end{equation}
So we have already seen in Theorem \ref{ap:th:3} that 
$y=z^{\lambda_0}\left[p_0+zp_1+p^2 y_2+ \cdots\right]$ is a 
solution. 

Again, we consider the formal power series \eqref{eq:power} 
with coefficients \eqref{eq:formal}. Differentiating 
\eqref{eq:formal-diff} with respect to $\lambda$, we have 
\begin{equation}\label{eq:sol1}
    \left[\frac{d}{dz}-A\right]\hat y(\lambda,y_0)=
         z^{\lambda-1}y_0 + z^{\lambda-1}\log z(\lambda I-R)y_0 ,
\end{equation}
where we put 
\[
 \hat y(\lambda,y_0)=\frac{\partial}{\partial \lambda}y(\lambda,y_0).
\]
In particular, $\hat y:=\hat y(\lambda_0,p_0)$ satisfies 
\[
    \left[\frac{d}{dz}-A\right]\hat y= z^{\lambda_0-1}p_0,
\]
because $(\lambda_0I-R)p_0=0$.  
In \eqref{eq:power}, if we put 
\[ 
      y(\lambda,p_0)=z^{\lambda}(p_0+zp_1+\cdots)
                        =:z^{\lambda}p(\lambda,z)
   \qquad\text{($p(\lambda,z)$ is holomorphic in $z$)}
\]
then 
\[
     \hat y(\lambda,p_0) 
         = \frac{\partial y}{\partial\lambda}(\lambda,z)
         = (\log z) z^{\lambda}p(\lambda,z)+
           z^{\lambda}\frac{\partial p}{\partial \lambda}(\lambda,z)
         = (\log z) y(\lambda,p_0)+
             z^{\lambda} q(\lambda,z),
\]
where $q(\lambda,z)(=\pd p/\pd \lambda)$ is a holomorphic 
function in $z$. 
Hence, it can be written as 
\begin{equation}\label{eq:yhat}
  \hat y= \hat y(\lambda_0,p_0)
        = (\log z) y+ z^{\lambda_0}q(\lambda_0,z). 
\end{equation}
On the other hand, if we set 
\begin{equation}\label{eq:defofw}
   w= y(\lambda_0, q_0),
\end{equation}
we have 
\[
  \left [\frac{d}{dz}-A \right] w=
z^{\lambda_0-1}(\lambda_0I-R)q_0=- z^{\lambda_0-1} p_0,
\]
from \eqref{eq:formal-diff} and \eqref{eq:eigen-vector}. 
Therefore, 
\[
   \tilde y:=\hat y + w
\]
is a solution to \eqref{eq:ode}. Moreover, $\tilde y$ is 
of the form $\tilde y=z^{\lambda_0}\left[q_0+zq_1+z^2 q_2+ \cdots\right]+
             (\log z)y$, by \eqref{eq:yhat} and \eqref{eq:defofw}.
\end{proof}

\normalsize

\end{document}